\documentclass{article}
\pdfoutput=1
\usepackage[numbers,sort&compress]{natbib}
\usepackage{graphicx, amssymb, amsmath}
\usepackage{amsthm}
\usepackage{mathrsfs}
\usepackage{color}
\usepackage{caption}
\usepackage{algorithm}
\usepackage{algorithmic}
\usepackage{upgreek}
\usepackage{subfigure}
\usepackage{color,soul}
\usepackage[colorlinks,
            linkcolor=blue,
            anchorcolor=blue,
            citecolor=blue]{hyperref}
\usepackage{subfigure}
\usepackage{setspace}
\sethlcolor{blue}
\usepackage{float}
\usepackage{natbib}
\usepackage{diagbox}
\setcitestyle{numbers,square}

\def\k{\kern .5em}
\def\er{\kern .2em}

\begin{document}

\date{}
\author{}
\newcommand{\be}{\begin{equation}}
\newcommand{\ee}{\end{equation}}
\newcommand{\ba}{\begin{array}}
\newcommand{\ea}{\end{array}}
\newcommand{\beas}{\begin{eqnarray*}}
\newcommand{\eeas}{\end{eqnarray*}}
\newcommand{\bea}{\begin{eqnarray}}
\newcommand{\eea}{\end{eqnarray}}
\newcommand{\ome}{\Omega}

\newtheorem{theorem}{Theorem}[section]
\newtheorem{lemma}{Lemma}[section]
\newtheorem{remark}{Remark}[section]
\newtheorem{proposition}{Proposition}[section]
\newtheorem{definition}{Definition}[section]
\newtheorem{corollary}{Corollary}[section]

\newtheorem{theo}{Theorem}[section]
\newtheorem{lemm}{Lemma}[section]
\newcommand{\blem}{\begin{lemma}}
\newcommand{\elem}{\end{lemma}}
\newcommand{\bthe}{\begin{theorem}}
\newcommand{\ethe}{\end{theorem}}
\newtheorem{prop}{Proposition}[section]
\newcommand{\bprop}{\begin{proposition}}
\newcommand{\eprop}{\end{proposition}}
\newtheorem{defi}{Definition}[section]
\newtheorem{coro}{Corollary}[section]
\newtheorem{algo}{Algorithm}[section]
\newtheorem{rema}{Remark}[section]
\newtheorem{property}{Property}[section]
\newtheorem{assu}{Assumption}[section]
\newtheorem{exam}{Example}[section]

\renewcommand{\theequation}{\arabic{section}.\arabic{equation}}
\renewcommand{\thetheorem}{\arabic{section}.\arabic{theorem}}
\renewcommand{\thelemma}{\arabic{section}.\arabic{lemma}}
\renewcommand{\theproposition}{\arabic{section}.\arabic{proposition}}
\renewcommand{\thedefinition}{\arabic{section}.\arabic{definition}}
\renewcommand{\thecorollary}{\arabic{section}.\arabic{corollary}}
\renewcommand{\thealgorithm}{\arabic{section}.\arabic{algorithm}}
\newcommand{\lan}{\langle}
\newcommand{\curl}{{\bf curl \;}}
\newcommand{\rot}{{\rm curl}}
\newcommand{\grad}{{\bf grad \;}}
\newcommand{\dvg}{{\rm div \,}}
\newcommand{\ran}{\rangle}
\newcommand{\bR}{\mbox{\bf R}}
\newcommand{\bRn}{{\bf R}^3}
\newcommand{\Coinf}{C_0^{\infty}}
\newcommand{\disp}{\displaystyle}
\newcommand{\ra}{\rightarrow}
\newcommand{\Ra}{\Rightarrow}
\newcommand{\ud}{u_{\delta}}
\newcommand{\Ed}{E_{\delta}}
\newcommand{\Hd}{H_{\delta}}
\newcommand\varep{\varepsilon}
\newcommand{\RNum}[1]{\uppercase\expandafter{\romannumeral #1\relax}}
\newcommand{\tabincell}[2]{\begin{tabular}{@{}#1@{}}#2\end{tabular}}
\title{Fast Algorithms for Finite Element nonlinear Discrete Systems to Solve the PNP Equations}

\author{*
	\and *
	\and *
}
\author{Yang Liu $^{1}$
	\and Shi Shu $^{2,*}$
	\and Ying Yang $^{3,*}$
}
\footnotetext[1]
{School of Mathematics and Computational Science, Xiangtan University,
	Xiangtan, 411105, P.R. China. E-mail: liuyang@smail.xtu.edu.cn
}

\footnotetext[2]
{$^{*}$\textbf{Corresponding author.} School of Mathematics and Computational Science, Xiangtan University,
	Xiangtan, 411105, P.R. China. E-mail: shushi@xtu.edu.cn }

\footnotetext[3]
{$^{*}$\textbf{Corresponding author.} School of  Mathematics and Computational Science, Guangxi Colleges and Universities Key Laboratory of Data Analysis and Computation, Guangxi Applied Mathematics Center (GUET), Guilin University of Electronic Technology, Guilin, 541004, Guangxi, P.R. China. E-mail: yangying@lsec.cc.ac.cn }

%
%

\maketitle
\linespread{1.5}
\noindent
{\bf Abstract}\quad The Poisson-Nernst-Planck (PNP) equations are one of the most effective model for describing electrostatic interactions and diffusion processes in ion solution
systems, and have been widely used in the numerical simulations of biological ion
channels, semiconductor devices, and nanopore systems. Due to the characteristics
of strong coupling, convection dominance, nonlinearity and multiscale, 
the classic Gummel iteration for the nonlinear discrete system of PNP equations
 converges slowly or even diverges.
We focus on fast algorithms of nonlinear discrete system for the general PNP equations, which have better adaptability, friendliness and efficiency than the Gummel iteration. First, a geometric full approximation storage (FAS) algorithm is proposed to improve the slow convergence speed of the Gummel iteration. Second, an algebraic FAS algorithm is designed, which does not require multi-level geometric information and is more suitable for practical computation compared with the geometric one.
Finally, improved algorithms based on the acceleration technique and adaptive method are proposed to solve 
the problems of excessive coarse grid iterations and insufficient adaptability to the size of computational domain in the algebraic FAS algorithm.
  The numerical experiments are shown for the geometric, algebraic FAS and improved algorithms respectively to illustrate the effiency of the algorithms.


\noindent
{\bf Keywords}: Poisson-Nernst-Planck equations, Gummel iteration, full approximation storage; adaptive and acceleration technology

\noindent
{\bf AMS(2000)subject classifications}\quad 65N30, 65N55

\section{Introduction}\label{sec1}
\noindent
{\color{black}The Poisson-Nernst-Planck (PNP) equations are a coupled nonlinear partial
	differential equations,
	which describe electrostatic interactions and diffusion processes in ion solution systems. They were first proposed by Nernst
	\textsuperscript{\cite{nernst1889elektromotorische}}
	and Planck\textsuperscript{\cite{planck1890ueber}} and have been widely used in the numerical simulations of biological ion channel\textsuperscript{\cite{tu2013parallel,eisenberg1998ionic,singer2009poisson}}, semiconductor devices\textsuperscript{\cite{markowich1985stationary,brezzi2005discretization,miller1999application}}, and nanopore systems\textsuperscript{\cite{xu2018time,daiguji2004ion}.}} Due to the strong coupling, convection dominance, nonlinearity, and multiscale characteristics of PNP practical problems, it is difficult to construct numerical methods with good approximation, stability and efficiency.
\vskip 0.3cm
Many numerical methods have been applied to solve the PNP equations, such as finite element (FE) \textsuperscript{
	\cite{lu2010poisson,xie2020effective,ying2021new,yang2020superconvergent}}, finite difference \textsuperscript{\cite{cardenas2000three,liu2014free,flavell2014conservative,he2016energy}}, and finite volume \textsuperscript{\cite{chainais2004finite,chainais2003convergence,chainais2003finite,bessemoulin2014study}} methods, etc.
Since there exist problems such as the poor stability or the bad approximation when using standard FE, finite difference or finite volume methods to discretize pratical PNP equations, some improved methods have emerged. The inverse average FE method was constructed for a class of steady-state PNP equations in nanopore systems in \cite{zhang2021inverse}, which effectively solves the problems of non physical pseudo oscillations caused by convection dominance. A new stable FE method called SUPG-IP was proposed in \cite{SUPG-IP} for the steady-state PNP equations in biological ion channel. Numerical experiments showed that the method has better performance than the standard~FE and~SUPG methods in positive conservation and robustness. The robustness and stability of the FE algorithm in channel numerical simulation were improved in \cite{tu2015stabilized} for the modified PNP equations considering ion size effects by using the SUPG method and pseudo residual free bubble function technique.
 In \cite{zhang2022class}, the~Slotboom variable was used to transform the continuity equation into an equivalent equation for a class of steady-state PNP equations in semiconductor, and then four FE methods with averaging techniques  were introduced  to discretize the nonlinear system.
Numerical experiments have shown that the new methods can ensure the conservation of the calculated terminal current and are as stable as the~SUPG method.

\vskip 0.3cm
The above nonlinear discrete systems are usually linearized by iterative methods such as~Gummel or~Newton, and then fast algorithms are designed for these linear systems. In \cite{mathur2009multigrid}, the Newton linearization was used for the adaptive finite volume discrete system of PNP equation on unstructured grid, and algebraic multigrid (AMG) method was applied to solve the linearized system.
 In \cite{bousquet2018newton},  a fast FGMRES algorithm was designed based on a block upper triangular preconditioner for the edge average finite element (EAFE) scheme of the Newton linearized PDEs of the PNP and Navier-Stokes coupled system, in which the~PNP subsystem was solved using the~PGMRES based on the~UA-AMG preconditioner.
  Xie and Lu \textsuperscript{\cite{xie2020effective}} used Slotboom transformation to transform the PNP equations with periodic boundary conditions into equivalent equations, and gave an acceleration method for the Gummel iteration, in which the nonlinear potential equation was solved by the modified Newton iterative method. Since the discrete system of the practical PNP problem has the characteristics of strong asymmetry, coupling between physical quantities, multi-scale, large scale and condition number, the above iterative algorithms may have slow convergence (even divergence) and low computational efficiency. There is another kind of method at present, whose basic idea is based on the two-grid method \textsuperscript{\cite{jin2006two}}, which essentially transforms the nonlinear discrete system on the fine grid to the nonlinear system on the coarse grid. Shen et. al. \textsuperscript{\cite{shen2020decoupling}} proposed a two-grid FE method for a class of time-dependent PNP equations. Based on the optimal $L^2$ norm error estimates of the derived discrete scheme, the~$H^1$ norm error estimates are presented.
Comparative experiments with Gummel iteration show that the new method has higher computational efficiency. Ying et. al. \textsuperscript{\cite{ying2021new}} aimed at the two-grid FE method of the~SUPG scheme for the steady-state~PNP equations in biological ion channel, in which the Krylov iterative method based on the block Gauss-Seidel preconditioner was designed for the Newton linearized system on the coarse grid, and it was proved that the condition number of the preconditioned system is independent of the mesh size under certain assumptions. It should be pointed out that although the above geometric two-grid method has higher computational efficiency than Gummel or Newton iterative method for some discrete systems of PNP equations, it may encounter many difficulties when the method is applied to solve complex practical PNP problems. For example, the approximate accuracy depends on the coarsening rate. Due to factors such as low smoothness of the solution, the coarsening rate cannot be too large, which affects the acceleration effect. In addition, complex meshes generated by professional/commercial softwares are usually difficult to meet the requirements of geometric two-grid method for nested meshes. Therefore, further research is needed on fast algorithms for PNP large-scale discrete systems.
\vskip 0.3cm
The FAS algorithm proposed by Brandt \textsuperscript{\cite{brandt1977multi}} is a nonlinear geometric multigrid method.  It has been applied to direct solutions of many nonlinear systems \textsuperscript{\cite{waisman2008heterogeneous,liu2017multigrid,feng2018mass,adams2010toward}}, the detailed descriptions of which refer to e.g.
\cite{brandt2011multigrid,briggs2000multigrid,brandt2011full}.
%
The FAS algorithm requires multi-level grid information, and at present coarse grids are mainly generated through aggregation methods \textsuperscript{
	\cite{alloush2022development,birken2019preconditioned}}. However, for complex domains and unstructured grids, aggregation algorithms lead to poor quality of coarse grids, which in turn affects the performance of the FAS algorithm\textsuperscript{
	\cite{alloush2022development}}. Therefore, it is important to study the algebraic FAS algorithm directly from the discretized nonlinear algebraic system of the PNP equations. {\color{black}In addition, due to the drawbacks of existing linearized methods for solving PNP discrete systems, such as insufficient adaptability to the size of domain
	and excessive iterations for solving corresponding nonlinear algebraic systems on coarse grid \textsuperscript{\cite{lu2007electrodiffusion,zhang2021inverse}}, further improvements on them will be worthwhile.}

 \vskip 0.3cm
  In this paper, we mainly study the fast algorithms of the finite element nonlinear discrete system to solve the PNP equations. First, for nonlinear discrete systems of general PNP equations, in order to overcome the drawbacks of slow convergence speed and weak adaptability to size of computational domain in the classical Gummel iteration, a geometric FAS algorithm is designed. Compared with the Gummel iteration, it not only has higher efficiency but also enhances the adaptability to convective intensity.
  Second, considering the limitations of the geometric FAS algorithm, such as the requirement of nested coarse grids, an algebraic FAS algorithm is presented. This algorithm does not require multi-level geometric information, which is more suitable for practical PNP problems. The numerical experiments including a practical problem in ion channel show that the algebraic FAS algorithm can not only improve the efficiency of decoupling external iterations, but also expand the domain of the problem that can be solved compared with the Gummel iteration, which indicates the effectiveness of this algorithm in solving the strong coupling and convection dominated PNP equations. Finally, some improved algorithms based on acceleration technique and adaptive method are designed to
  deal with the problems of excessive coarse grid iterations and insufficient adaptability to the size of computational domain in the algebraic FAS algorithm when the convective intensity is high. Numerical experiments show that the improved algorithms converge faster and significantly enhances the adaptability of convective intensity.
\vskip 0.3cm
The rest of paper is organized as follows. In Section 2, we introduce some standard notations and PNP equations. The continuous and discrete forms are also given in this section. In Section 3, we show the geometric and algebraic FAS algorithms. The improved algrithms based on the acceleration technique and adaptive method are also presented in this section. The conclusion is presented in Section 4.

\setcounter{lemma}{0}
\setcounter{theorem}{0}
\setcounter{corollary}{0}
\setcounter{equation}{0}
\setcounter{remark}{0}
\section{ The continuous and discrete problems}
Let $\Omega  \subset \mathbb{R}^{d}$ $(d=2,3)$ be a bounded Lipschitz domain. We clarify the standard notations for Sobolev spaces $W^{s,p}(\Omega)$ and
their associated norms and seminorms. For $p=2$, we denote
$H^{s}(\Omega)=W^{s,2}(\Omega)$ and $H_0^1(\Omega)=\{v\in H^{1}(\Omega):v|_{\partial\Omega}=0\}$. For simplicity, let $\|\cdot\|_{s}=\|\cdot\|_{W^{s,2}(\Omega)}$ and $\|\cdot\|=\|\cdot\|_{L^2(\Omega)}$ and $(\cdot,\cdot)$ denote the standard $L^{2}$-inner product.
\vspace{3mm}
\subsection{PNP equations}

 We consider the following time-dependent PNP equations (cf. \cite{shen2020decoupling})
 \begin{equation}\label{pnp equation}
 \left\{\begin{array}{lr}
 -\Delta\phi-\sum\limits^{2}_{i=1}q^{i}p^{i}=f, & \text { in } \Omega,~\text{for}~t\in(0, T], \vspace{1mm}\\
 {\partial_t p^{i}}-\nabla\cdot(\nabla{p^{i}}+ q^{i}p^{i}\nabla\phi)=F^i, & \text { in } \Omega,~\text{for}~t\in(0, T], i=1,2,\vspace{1mm}
 \end{array}\right.
 \end{equation}
 where $p^i,~i=1,2$ is the concentration of the $i$-th ionic species with charge $q^i$, $p_t^{i}=\frac{\partial p^i}{\partial t}$, $\phi$ denotes the electrostatic potential, 
$f$ and $F^{i}$ are the reaction source terms,
 The homogeneous Dirichlet boundary conditions are employed
 \begin{equation}\label{boundary}
 \left\{
 \begin{array}{c}
 \;\phi=0,\;\text{on}\;\partial\Omega,~\text{for}~t\in(0, T],\vspace{0.5mm}\\
 \;p^i=0,\;\text{on}\;\partial\Omega,~\text{for}~t\in(0, T].
 \end{array}
 \right.
 \end{equation}

 The corresponding steady-state PNP equations are as follows (cf. \cite{liu2021virtual})
 \begin{equation}\label{statPNP}
 \left\{\begin{array}{lr}
 -\Delta\phi-\sum\limits^{2}_{i=1}q^{i}p^{i}=f, & \text { in } \Omega,~\text{for}~t\in(0, T], \vspace{1mm}\\
 -\nabla\cdot(\nabla{p^{i}}+ q^{i}p^{i}\nabla\phi)=F^i, & \text { in } \Omega,~\text{for}~t\in(0, T], i=1,2,\vspace{1mm}
 \end{array}\right.
 \end{equation}

 The weak formulations of (\ref{pnp equation})-(\ref{boundary})  and (\ref{statPNP})-(\ref{boundary}) are that:
 find $p^i\in L^2(0,T;H_0^1(\Omega))$ $\cap L^{\infty}(0,T;L^{\infty}(\Omega))$, $i=1,2$, and $\phi(t)\in H_0^1(\Omega)$ such that
 \begin{eqnarray}
 \label{weak1con}
 (\nabla\phi,\nabla w)-\sum\limits^{2}_{i=1}q^{i}(p^{i},w)=(f,w),~~\forall w\in H_0^1(\Omega), \\
 \label{weak2con}
 (\partial_tp^i,v)+(\nabla{p^{i}},\nabla v)+ (q^{i}p^{i}\nabla\phi,
 \nabla v)=(F^i,v), ~~\forall v\in H_0^1(\Omega), ~i=1,2.
 \end{eqnarray}
 and find $p^i  \in H_0^1 (\Omega ),~i = 1,2$ and $ \phi  \in H_0^1 (\Omega )$ such that
 \begin{eqnarray}
 \label{statweak1con}
 (\nabla\phi,\nabla w)-\sum\limits^{2}_{i=1}q^{i}(p^{i},w)=(f,w),~~\forall w\in H_0^1(\Omega), \\
 \label{statweak2con}
 (\nabla{p^{i}},\nabla v)+ (q^{i}p^{i}\nabla\phi,
 \nabla v)=(F^i,v), ~~\forall v\in H_0^1(\Omega), ~i=1,2,
 \end{eqnarray}
 respectively.

  {The existence and uniqueness of the solutious to (\ref{weak1con})-(\ref{weak2con}) and (\ref{statweak1con})-(\ref{statweak2con}) have been presented in \cite{gajewski1986basic} for $F^i=r(p^1,p^2)(1-p^1p^2)$.
  	Here $p^1,~p^2$ represent the densities of mobile holes and electrons respectively in a semiconductor device. The function $r: R_+^2\rightarrow R_+$ are required to be Lipschitzian and $p^1,~p^2\in W^{0,\infty}(\Omega)$.

 \subsection{The discrete scheme}

 Suppose $\mathcal{T}_h = \{ K \}$ is a  partition of $\Omega$, where $K$ is the element and $h=\max\limits_{K\in\mathcal{T}_h} \{\mbox{diam}~K\}$. We define the first order finite element subspace of $H_0^1(\Omega)$ as follows
 \begin{eqnarray}
 V_h =\{v \in H^1(\Omega): v|_{\partial\Omega}=0~\mbox{and}~v|_K \in \mathcal{P}_1(K), ~~\forall K\in \mathcal{T}_h\},
 \end{eqnarray}
 where $\mathcal{P}_1(K)$
 is the set of linear polynomials on the element $K$.

  The semi-discretization to \eqref{weak1con}-\eqref{weak2con} is defined as follows: find $p^{i}_h,~i=1,2$ and $\phi_h\in V_h$ such that
 \begin{eqnarray}
 \label{semiweak2}
 (\nabla\phi_h,\nabla w_h)-\sum\limits^{2}_{i=1}q^{i}(p^{i}_h,w_h)=(f,w_h),~~\forall w_h\in V_h, \\
 \label{semiweak1}
(\partial_tp^i_h,v_h)+(\nabla{p^{i}_h},\nabla v_h)+ (q^{i}p^{i}_h\nabla\phi_h,
 \nabla v_h)=(F^i,v_h), ~~\forall v_h\in V_h,~i=1,2.
 \end{eqnarray}

 In order to present the full discrete approximation to  (\ref{pnp equation})-\eqref{boundary}, let $0=t^0<t^1<\cdots<t^N=T$ denote the partition of $[0,T]$ into subintervals $(t^{n-1},t^n)$ such that $\tau=\max\{t^n-t^{n-1},n=1,2,\cdots N\}$. For any function $u$, denote by
 \begin{align*}
 u^n(x)=u(x,t^n)
 \end{align*}
 and
 $$D_\tau u^{n+1} =\frac{u^{n+1}-u^{n}} { \tau} , ~\mbox{for}~ n=0,1,2,\cdots,N-1.$$
The backward Euler full discretization based on FE method for  (\ref{pnp equation})-\eqref{boundary} is: Find
$p_h^{i,n+1}\in V_h,~i=1,2$ and $\phi_h^{n+1}\in V_h$, such that
{\small
\begin{align} \label{fulldis2-FEM}
 (\nabla\phi_h^{n+1},\nabla w_h)-\sum\limits_{i = 1}^2 q^i(p_h^{i,n+1},w_h)
 &= (f^{n+1 },w_h),  \forall w_h\in V_h,\\\label{fulldis-FEM}
 (D_\tau p_h^{i,n+1},v_h)+(\nabla p_h^{i,n+1},\nabla v_h)+(q^ip_h^{i,n+1}\nabla\phi^{n+1}_h,\nabla v_h)
 &= (F^{i,n+1},v_h), \forall v_h\in V_h,
 \end{align}}
where $f^{n+1 }=f (t^{n+1},\cdot)  $, $F^{i,n+1}=F^i(t^{n+1},\cdot),~i=1,2$.

There exist some convection dominated PNP equations in semiconductor area, for which the standard FE method does not work well.
The EAFE method has shown advantages in solving convection dominated equations. In order to present the EAFE scheme for PNP equations, suppose $K$ denotes the element with the number $k$, and $E$ is the edge with the endpoints $x_{k_\nu}$ and $x_{k_\mu}$, where ${k_\nu}$ is the whole number corresponding to the local number $\nu$. Let $\tau_E= x_{k_\nu}-x_{k_\mu},~{k_\nu} <{k_\mu} $ for any $E$ in $\mathcal T_h$.


 The EAFE fully discrete scheme for PNP equations is: find
 $p^{i,n+1}_h\in  V_h ,~i=1,2$  and $\phi_h^{n+1}\in V_h$, satisfying
 {\small\begin{align}
 \label{fulldis2-EAFE}
 (\nabla\phi_h^{n+1},\nabla w_h)&-\sum\limits_{i = 1}^2 q^i(p_h^{i,n+1},w_h)
=   (f^{n+1 },w_h),~\forall w_h\in V_h ,\\
(D_\tau p_h^{i,n+1},v_h) &+ \sum\limits_{K\in \mathcal{T}_{h}}[\sum\limits_{E\subset K}\omega_E^K\tilde{\alpha}^{K,i}_{E}(\phi_h^{n+1})\delta_E(e^{q^i\phi_h^{n+1}}p_h^{i,n+1})\delta_Ev_{h}]
= \sum\limits_{K\in \mathcal{T}_{h}}(F^{i,n+1},v_h)_K,  ~\forall v_h\in V_h ,  \label{fulldis-EAFE}
 \end{align}}
 where
{\small\begin{equation}\label{omega}
 	\tilde{\alpha}^{K,i}_{E }(\phi_h^{i,n+1} )=\left(\frac{1}{|\tau_{E }|}\int_{E } e^{q^i\phi^{n+1} _h}\mathrm{d}s\right)^{-1},~
 	\omega_{E }^K  :=\omega_{E_{\nu \mu }}^K  =-(\nabla \psi_{k_\mu},\nabla \psi_{k_\nu})_K ,~\delta_{E }( u)=  u_{k_\nu}- u_{k_\mu}, ~ u=e^{q^i\phi_h }p_h^{i}.
\end{equation}}

To present the corresponding discrete algebraic system, we define some notations. Let
\begin{eqnarray}\label{Psi}
\Psi=(\psi_{1},\ldots,\psi_{{n_h}}),~{n_h}=\dim(V_h )
\end{eqnarray}
be a set of basis function vectors of $V_h$.
Note that $\forall v_h\in V_h$, we have
\begin{equation}\label{phihEi}
v_h|_K=\sum_{m=1}^4V_{k_m}\lambda^K_m= \lambda^K  V_K,
\end{equation}
where $ \lambda^K =(\lambda^K_1,\lambda^K_2,\lambda^K_3,\lambda^K_4)$ and $ V_K=(V_{k_1},V_{k_2},V_{k_3},V_{k_4})^T $ are the volume coordinate vector and degree of freedom vector of element K, respectively.
Assume $n_h$ dimensional vectors
\begin{equation}\label{def-PhiG3}
\Phi^{n+1}=(\Phi^{n+1}_1,\ldots,\Phi^{n+1}_{n_h})^T,\quad G_\Phi^{n+1}=((f^{n+1 },\psi_1),\ldots,(f^{n+1 },\psi_{n_h}))^T,
\end{equation}
and $2n_h$ dimensional vectors
\begin{equation}\label{def-Pn}
P^\mu=
\begin{pmatrix}
P^{1,\mu}\\
P^{2,\mu}
\end{pmatrix}
,~~
P^{i,\mu}=(p^{i,\mu}_1,\ldots,p^{i,\mu}_{n_h})^T,~\mu=n,n+1,~ i=1,2,
\end{equation}
{\small
	\begin{equation}\label{def-PG}
	G^{n+1}=
	\begin{pmatrix}
	G^{1,n+1}\\
	G^{2,n+1}
	\end{pmatrix},~ G^{i ,n+1}=(g^{i ,n+1}_1,\ldots,g^{i ,n+1}_{n_h})=((F^{i ,n+1},\psi_1),\ldots,(F^{i,n+1},\psi_{n_h}))^T,~ i=1,2,
	\end{equation}}
$n_h\times 2n_h$ lumped mass matrix
\begin{equation}  \label{C4-M-def}
\bar M =-(q^1M,q^2 M), ~ M=\frac{1}{4}\textrm{diag}(|\Omega_1|,\cdots ,|\Omega_{n_h}|), ~ \Omega_k=supp(\psi_k),~k=1,\ldots,{n_h},
\end{equation}
and $n_h\times n_h$ stiffness matrix
\begin{equation}  \label{C4-AL-def}
A_L=(  \nabla\Psi^T,\nabla \Psi).
\end{equation}
%
%

It is easy to know the nonlinear algebraic equation for the EAFE discrete system ~\eqref{fulldis2-EAFE}-\eqref{fulldis-EAFE} is as  follows:
\begin{equation}\label{1.1-EAFE}
A_1(U_1)U_1=F_1.
\end{equation}
Here
	\begin{equation}\label{def-A1U1F1}
	A_1(U_1):=\left(
	\begin{array}{cc}
	A_{1}^\Phi& \bar M_1\\
	\bf 0    &  A_1^{P}(\Phi_1) \\
	\end{array}
	\right)
	,~U_1=\left(
	\begin{array}{c}
	\Phi_1\\P_1
	\end{array}
	\right),~F_1=\left(
	\begin{array}{c}
	F_1^\Phi \\
	F_1^{P}
	\end{array}
	\right),
	\end{equation}
where
\begin{equation}\label{def-A1F1}
A_{1}^\Phi=A_L,~ \bar M_1=\begin{pmatrix}
-q^1M_1& -q^2 M_1
\end{pmatrix},~ M_1=  M,~A_1^{P}(\Phi_1):= {\bf M}+ \tau \bar A (\Phi^{n+1} ),
\end{equation}
\begin{equation}\label{def-U1}
\Phi_1=\Phi^{ n+1 },~P_1=P^{ n+1 },~
F ^\Phi=G_\Phi^{n+1},~F ^{P}=F^{n }.
\end{equation}
matrix
${\bf M} =diag(M, M)$, $ \bar M$ and $A_L$ are defined by \eqref{C4-M-def} and \eqref{C4-AL-def}, respectively, $2n_h\times 2n_h$ stiffness matrix
\begin{equation}  \label{C4-A-def}
\bar A (\Phi^{n+1})=diag(\bar A^1(\Phi^{n+1}),\bar A^2(\Phi^{n+1})),
\end{equation}
the general element of the element stiffness matrix ~$\bar A^{i,K}(\Phi^{n+1})$  on tetrahedron element $K$ for $\bar A^i(\Phi^{n+1})$ is as follows:
{
	\begin{equation}\label{aij}
	\bar a_{\nu \mu }^{i,K}(\Phi^{n+1})=\left\{\begin{array}{ll}
	& -\omega^K_{E_{\nu \mu }}\tilde{\alpha} _{E_{\nu \mu }}^{K,i}(\Phi^{n+1})e^{q^i \lambda^K(q_{\mu })  \Phi ^{n+1}_K} , {\nu >\mu }, \\
	& -\omega^K_{E_{\mu \nu }}\tilde{\alpha} _{E_{\mu \nu }}^{K,i}(\Phi^{n+1})e^{q^i \lambda^K(q_{\mu })  \Phi ^{n+1}_K} , {\nu <\mu }, \\
	&\sum\limits_{ k> \mu}\omega^K_{E_{\mu  k}}\tilde{\alpha} _{E_{\mu k}}^{K,i}(\Phi^{n+1})e^{q^i \lambda^K(q_{\mu })  \Phi ^{n+1}_K}
	\\&+ \sum\limits_{ k<\mu} \omega^K_{E_{k\mu}}\tilde{\alpha} _{E_{k\nu_2}}^{K,i}(\Phi^{n+1})e^{q^i \lambda^K(q_{\mu })  \Phi ^{n+1}_K} ,  {\nu =k},
	\end{array}\right. ~\nu ,\mu=1,2,3,4,~i=1,2,
	\end{equation}}
where~$\omega^K_{E_{\nu \mu }}$ is defined by~\eqref{omega}, $q_1,q_2,q_3,q_4$ are four vertices of element~$K$, and
\begin{equation}\label{[1](3.10)}
\tilde{\alpha}^{K,i}_{E_{\nu \mu } }(\Phi^{n+1} )=\left[\frac{1}{|\tau_{E }|}\int_{E } e^{q^i \lambda^K  \Phi ^{n+1}_K}\mathrm{d}s\right]^{-1}.
\end{equation}

 \setcounter{lemma}{0}
 \setcounter{theorem}{0}
 \setcounter{corollary}{0}
 \setcounter{equation}{0}
  \section{Fast algorithms for PNP equations }
 In this section, first a geometric FAS algorithm is proposed, which can overcome the slow convergence speed of the classical Gummel algorithm for the nonlinear discrete systems of general PNP equations. Then due to the limitation of the geometric FAS algorithm, such as the requirement for nested coarse grids, an algebraic FAS algorithm is designed. However,  the algebraic FAS algorithm still has drawbacks such as too many coarse grid iterations when encountering the problem of strong convection dominance. Hence, we  design some improved algorithms based on the acceleration technique and adaptive method.

 \subsection{Geometric and algebraic FAS algorithms}
In this subsection, we first present the geometric FAS algorithm, then design algebraic FAS algorithm. Numerical examples show the FAS algorithm improves the efficiency of the  decoupling iterations and expands the computational domain of PNP equations.
}

%
%

\subsubsection{Geometric FAS algorithm}
The geometric FAS algorithm for PNP equations is as shown in Algorithm \ref{g-FAS}, which can be constructed by following \cite{brandt1977multi,bueler2021full} and using the Gummel iteration for PNP equations \textsuperscript{\cite{liu2021virtual}}.
\begin{algorithm}[!h]
	\caption{Geometric FAS Algorithm}\label{g-FAS}
	\begin{description}
		\item[Setup:]~
		\begin{description}
			\item[S1.] Generate two nested grids   $\mathcal{T}_{h_l},~l=1,2,~h_1<h_2$.  Let $V_{h_l}$ be the corresponding linear FE function space,  $n_l=\dim(V_{h_l})$;   Take prolongation operator $P_2^1$ as the linear FE interpolantion function from $V_{h_2}$ to $V_{h_1}$. 			
			
			\item[S2.]Generate restriction operator $R_1^{2,u}: V_{h_1}\rightarrow V_{h_2}$.
			The general element of the corresponding matrix $R_1^{2,u}\in R^{n_{2}\times n_1}$ is defined by
			\begin{equation}\label{Ru2}
			R_1^{2,u}(i,j)=\left\{
			\begin{array}{ll}
			\frac{R_1^{2}(i,j)}{ \sum\limits_{j=1}^{n_{1}} R_1^{2}(i,j)}, & \hbox{$R_1^{2}(i,j) \ne 0$;} \\
			0, & \hbox{else,}
			\end{array}
			\right.,
			\end{equation}
			where~$ R_1^{2}=(P_2^1)^T$ (see Remark \ref{remark4.1} for the reason why two sets of restriction operators are needed.)
			
			\item[S3.] Obtain the submatrix $A^\Phi_l$ and $\bar M_l$ of the stiffness matrix $A_{l}(\cdot)$ corresponding to  \eqref{fulldis2-EAFE}-\eqref{fulldis-EAFE} on grids $\mathcal{T}_{h_l},l=1,2$, and load vector $F_{l}$, and give the intial value $U_1^{(0)}$, the stopping tolerance  $\epsilon$ and $\epsilon_c$, and the maximum number $M_{iter}$ of iterations on the coarse and fine grids.
		\end{description}
		
		\item[Solve:]
		
		For any $l\ge0$,  The FAS iteration from $U_{1}^{(l)}\rightarrow U_{1}^{(l+1)}$ is as follows.
		\begin{description}
			\item[S1.]{\bf Pre-smoothing}: Take $U_{1}^{(l)}$ as the initial value, and then perform $\nu_1 \ge 0$ iterations by using Gummel iteration.
			\item[S2.] {\bf Solving coarse grid system}:
			Solve the following nonlinear system on the coarse grid using Gummel iteration 
			\begin{equation}\label{ce-g}
			A _2\left(  U_2\right)  U_2= \tilde{\tau}_2.
			\end{equation}
			Here
			\begin{equation}\label{r1}
			\tilde{\tau}_2:=A_2\left(y_2\right)y_2 +\vec R_1^{2}r_1
			,~
			y_2 = \vec R_1^{2,u}   U_1^{(  f s)}, ~r_1=F_1-A_1\left(U_1^{(fs)}\right)  U_1^{(fs)},
			\end{equation}
			where
			\begin{equation}
			\vec R_1^{2} =diag(    R_1^{2}, R_1^{2}, R_1^{2}) ,~
			\vec R_1^{2,u} =diag(  R_1^{2,u}, R_1^{2,u}, R_1^{2,u}) .
			\end{equation}

			\item[S3.] {\bf Correction}:
			\begin{equation}
			U_1^{ ( bs)}=  U_1^{ (  fs )} + (\vec R_1^{2})^T(  U_2- y_2 ).
			\end{equation}
			\item[S4.] {\bf Post-smoothing}: Choosing $  U_{1}^{( bs)} $ as the initial value, perform $\nu_2\ge 0$ iterations for \eqref{1.1-EAFE}  by using Gummel iteration and get $U_{1}^{(l+1)}$.
		\end{description}
	\end{description}
\end{algorithm}

\begin{remark} \label{remark4.1}
	The matrix $R_1^{2,u}$ is the transfer matrix between degrees of freedom from $V_{h_1}$  to $V_{h_2}$. The matrix $R_1^{2}$ is the transfer matrix between the basis funcitons from $V_{h_1}$  to $V_{h_2}$. Since the row sum of $R_1^{2}$ may not be $1$, e.g. it may not satisfies the constant approximation, it can not be the trasfer operator between degrees of freedom and it needs to generate $R_1^{2,u}$ to transfer degrees of freedom which satisfies constant approximation.
\end{remark}

For given tolerance $\epsilon$, the exit criterion for Algorithm \ref{g-FAS} is $l+1\ge$ $M_{iter}$ or
\begin{align}\label{tol-fas}
\|F_1-A_1 (  U_1^{ ( l+1 )} )  U_1^{ ( l+1 )}  \|_{R^{3n_1}}\leq\epsilon.
\end{align}


If take $U_1^{(0)}=\bf 0$, $\nu_1=0$ in the geometric FAS algorthm (Algorithm \ref{g-FAS}), then from  $r_1=F_1=((f_3,\Psi_1),(f_1,\Psi_1),(f_2,\Psi_1))^T$  and $\Psi_2= R_1^{2}\Psi_1$, where $ \Psi_i,~i=1,2$ is
the FE basis function vector on grid $\mathcal T_{h_i}$, we have
$$\tilde \tau_2=\vec R_1^{2}r_1=((f_3,R_1^{2}\Psi_1),(f_1,R_1^{2}\Psi_1),(f_2,R_1^{2}\Psi_1))^T
=((f_3,\Psi_2),(f_1,\Psi_2),(f_2,\Psi_2))^T=F_2,$$
where $F_2$ is right hand of the coarse grid equation in the two-grid (TG) algorithm provided in  \cite{shen2020decoupling,ying2021new}. Hence, the TG algorithm can be viewed as a degradation of Algorithm
\ref{g-FAS} in the case when $U_1^{(0)}=\bf 0$, $\nu_1=0$, $\nu_2=1$.

Next we present a numerical example to observe the efficiency of geometric FAS algorithm and compare it with the Gummer iteration \cite{yang2020superconvergent} and TG algorithm \cite{shen2020decoupling,ying2021new}. First, we study the approximation effect of the Gummel iteration  and TG algorithm
under different convection intensity. Then we present the results of the geometric FAS algorithm and compare it with the Gummel and TG algorithms. \\

\begin{exam}\label{exam1}
	Consider the following steady-state PNP equations
	\begin{equation}\label{cont-new-model-steady}
	\left\{
	\begin{array}{llll}
	&-\Delta u- ( p-n)=f, & & \text { in } \Omega, \\
	&-\nabla \cdot\left(\nabla p +c_{\lambda}  p  \nabla u\right) =F^p,  & & \text { in } \Omega, \\
	&-\nabla \cdot\left(\nabla n-c_{\lambda} n \nabla u\right) =F^n,  & & \text { in } \Omega,\end{array}
	\right.
	\end{equation}
	with boundary conditions
	\begin{equation}\label{boundary-steady}
	\left\{
	\begin{array}{l}
	u=g_u,\text{on}\;\partial\Omega,\vspace{0.5mm}\\
	p =g_{p },\;\text{on}\;\partial\Omega,~i=1,2,\\
	n=g_{n},\;\text{on}\;\partial\Omega,~i=1,2.
	\end{array}
	\right.
	\end{equation}
Here $\Omega=[-0.5,0.5]^3$, $c_{\lambda}=0.179L^2=1.79$, the right hand function and the boundary value are given by the following exact solution
\begin{equation}\label{ue-beachmark}
\left\{
\begin{array}{ll}
u = \cos(\pi x)\cos(\pi y)\cos(\pi z),&\\
p  =3\pi^2(1+\frac{1}{2}\cos(\pi x)\cos(\pi y)\cos(\pi z)), &  \\
n=3\pi^2(1-\frac{1}{2}\cos(\pi x)\cos(\pi y)\cos(\pi z)). &
\end{array}
\right.
\end{equation}
	\end{exam}
The meshes used are uniform tetrahedral meshes, the tolerance $\epsilon=1E-6$, the maximum iteration number $M_{iter}$=1000, linear systems are solved by  the direct method solver Pardiso. For similicity, denote by
$e_{\zeta}=u_h-u$ and $\zeta=u,p ,n$.

Due to the fact that $c_\lambda=0.179L^2$ in \eqref{cont-new-model-steady} represents the convection coefficient, we show the errors of the Gummel iteration (see Tables \ref{Gummel-8}-\ref{Gummel-32}) and the number of iteration  (see Table \ref{gummel-iter}) with different $L$ to observe the impact of convection intensity on the Gummel iteration. It is shown on
Tables \ref{Gummel-8}-\ref{Gummel-32} show the $L^2$ and $H^1$ norm errors with different $L$, which indicates that the Gummel iteration is effective to solve PNP equations on uniform tetrahedral meshes when $L$ is not large.
 Howevever, if $L$ gets larger, which means strong convection dominance appears, the Gummel iterations converges very slowly even diverges, see Table  \ref{gummel-iter}.
\begin{table} [H]
	\caption{ Error of Gummel iteration when $h=\frac{1}{8}$  }\label{Gummel-8}
	\centering
	\begin{tabular}{|c|c|c|c|c|c|c|c|c| }
		\hline
		$L$ &$\|e_{ u }\|_0$   &  $\|e_{ p  }\|_0$ &  $\|e_{ n }\|_0$ & $\|e_{ u }\|_1$   & $\|e_{ p  }\|_1$ & $\|e_{ n }\|_1$   \\
		\hline
		$1$  & 2.32E-02&	4.00E-01&	4.06E-01&	4.80E-01&	7.10 &	7.10
		\\
		\hline
		$\sqrt{2.6}$ &2.43E-02	&4.08E-01&	4.32E-01&	4.80E-01&	7.11 	&7.11
		\\
		\hline
		$\sqrt{2.7}$  & 2.43E-02&	4.08E-01	&4.33E-01&	4.80E-01&	7.11 &	7.11
		\\
		\hline
	\end{tabular}
\end{table}
\begin{table} [H]
	\caption{ Error of Gummel iteration when $h=\frac{1}{16}$ }\label{Gummel-16}
	\centering
	\begin{tabular}{|c|c|c|c|c|c|c|c|c| }
		\hline
		$L$ & $\|e_{ u }\|_0$   &  $\|e_{ p  }\|_0$ &  $\|e_{ n }\|_0$ & $\|e_{ u }\|_1$   & $\|e_{ p  }\|_1$ & $\|e_{ n }\|_1$  \\
		\hline
		$1$  & 6.04E-03&	1.03E-01&	1.05E-01	&2.43E-01&	3.60 &	3.60
		\\
		\hline
		$\sqrt{2.6}$ & 6.12E-03	&1.08E-01&	1.14E-01&	2.43E-01	& 3.60 &	3.60
		\\
		\hline
		$\sqrt{2.7}$  &6.36E-03	&1.04E-01	&1.11E-01&	2.43E-01&	3.60 &	3.60
		\\
		\hline
	\end{tabular}
\end{table}
\begin{table} [H]
	\caption{ Error of Gummel iteration when $h=\frac{1}{32}$  }\label{Gummel-32}
	\centering
	\begin{tabular}{|c|c|c|c|c|c|c|c|c| }
		\hline
		$L$ &  $\|e_{ u }\|_0$   &  $\|e_{ p  }\|_0$ &  $\|e_{ n }\|_0$ & $\|e_{ u }\|_1$   & $\|e_{ p  }\|_1$ & $\|e_{ n }\|_1$   \\
		\hline
		$1$  & 2.16E-03&	2.25E-02&	2.34E-02&	1.22E-01&	1.80 &	1.80
		\\
		\hline
		$\sqrt{2.6}$ &2.46E-03&	1.65E-02&	1.72E-02&	1.22E-01&	1.80 &	1.80
		\\
		\hline
		$\sqrt{2.7}$  & 8.65E-04&	3.98E-02	&4.16E-02&	1.22E-01&	1.81 &	1.81
		\\
		\hline
	\end{tabular}
\end{table}

\begin{table} [H]
	\caption{The number of iterations of Gummel iteration}\label{gummel-iter}
	\centering
	\begin{tabular}{|c|c|c|c|c|c|c|c|c|c| }
		\hline
		\diagbox{$h$}{$L$}  & 1  & $\sqrt{2.6}$ &  $\sqrt{2.7 }$ &  $\sqrt{2.8}$  \\
		\hline
		$\frac{1}{8}$  & 11	&175&557&*
		\\
		\hline
		$\frac{1}{16}$ &9 &	118 &273&*
		\\
		\hline
		$\frac{1}{32}$ & 7	& 85 & 188  &*
		\\ \hline
	\end{tabular}
\end{table}
\noindent
where $*$ means the algorithm diverges.

Next, we compare the result of Gummel iteration with that of TG algorithm. Tables \ref{Gummel-TG-1}-\ref{Gummel-TG-2} display the results of the two algorithms on uniform tetrahedral meshes, where the coarsening rates are 1:8 and 1:64 respectively when the sizes of the coarse and fine grid are $\frac{1}{16}(\frac{1}{8})$ and  $\frac{1}{16}(\frac{1}{4})$, respectively.

\begin{table} [H]
	\caption{ Errors for Gummel and TG algorithms with $L=1$ }
	\label{Gummel-TG-1}
	\centering
	\begin{tabular}{|c|c|c|c|c|c|c|c|c|c| }
		\hline
		$h(H)$  &  $\|e_{ u }\|_0$   &  $\|e_{ p }\|_0$ &  $\|e_{  n }\|_0$ & $\|e_{ u }\|_1$   & $\|e_{ p  }\|_1$ & $\|e_{ n }\|_1$ \\
		\hline
		Gummel-$\frac{1}{16}$  &1.03E-02&7.49E-02&8.03E-02&2.44E-01&3.60&3.60
		\\
		\hline
		TG-$\frac{1}{16}(\frac{1}{8})$  & 2.34E-02		&4.59E-02	&4.47E-02&2.63E-01 &		3.64&		3.62
		\\
		\hline
		TG-$\frac{1}{16}(\frac{1}{4})$  & 1.16E-02 &	2.97E-01&	2.16E-01&4.43E-01	&	4.18 &3.97
		\\
		\hline
	\end{tabular}
\end{table}

\begin{table} [H]
	\caption{  Errors for Gummel and TG algorithms with $L=\sqrt{2.7}$   }
	\label{Gummel-TG-2}
	\centering
	\begin{tabular}{|c|c|c|c|c|c|c|c|c|c| }
		\hline
		$h(H)$  &  $\|e_{ u }\|_0$   &  $\|e_{ p }\|_0$ &  $\|e_{  n }\|_0$ & $\|e_{ u }\|_1$   & $\|e_{ p  }\|_1$ & $\|e_{ n }\|_1$  \\
		\hline
		Gummel-$\frac{1}{16}$  &9.09E-03&	5.47E-02&	7.49E-02&	2.43E-01&	3.61&	3.60
		\\
		\hline
		TG-$\frac{1}{16}(\frac{1}{8})$  &1.79E-02	&9.67E-02&	8.40E-02	&2.52E-01	&3.73&	3.71
		\\
		\hline
		TG-$\frac{1}{16}(\frac{1}{4})$  &5.34E-02&	6.13E-01&	5.69E-01&	3.67E-01&	5.40&	5.20
		\\
		\hline
	\end{tabular}
\end{table}

It is seen from Tables \ref{Gummel-TG-1}-\ref{Gummel-TG-2} that the $H^1$ and $L^2$ norm errors of the two algorithms are similar when the size $L$ of the computational domain is small. However, the errors of the TG algorithm is larger than that of the Gummul algorithm when $L$ is larger. For example, the errors $\| e_{p}\|_0$ of the TG algorithm with coarsening rates 1:8 and 1:64 are about 2 and 10 times that of the Gummel algorithm, respectively when $L=\sqrt{2.7}$.

Next, we investigate the adaptability of the geometric FAS algorithm to the size $L$. Tables \ref{16-fas-g}-\ref{32-fas-g} show the number of iterations and errors of algorithm \ref{g-FAS} on uniform tetrahedral meshes with a coarsening rate of 1:8, i.e., a mesh size of $\frac{1}{16}(\frac{1}{8})$ . In the computation, the number of pre-smoothing and post-smoothing times is $\nu_1=\nu_2=1$, stop tolerances $\epsilon=1E-6$ and $\epsilon_c=1E-7$, iter and iter-c respectively represent the number of iterations and the average number of coarse grid iterations of the geometric FAS algorithm.

\begin{table} [H]\small
	\caption{Errors of geometric FAS algorithm when $h(H)=\frac{1}{16}(\frac{1}{8})$  }\label{16-fas-g}
	\centering
	\begin{tabular}{|c|c|c|c|c|c|c|c|c| }
		\hline
		$L$ &    iter(iter-c) &    $\|e_{ u }\|_0$   &  $\|e_{ p }\|_0$ &  $\|e_{  n }\|_0$ & $\|e_{ u }\|_1$   & $\|e_{ p  }\|_1$ & $\|e_{ n }\|_1$     \\
		\hline
		${1}$ &3(7)&1.03E-02	&7.49E-02	&8.03E-02&	2.44E-01&	3.60&	3.60
		\\
		\hline
		$\sqrt{1.5}$ & 3(10)&  9.80E-03	&6.96E-02&	7.56E-02&	2.44E-01&	3.60&	3.60
		\\
		\hline
		$\sqrt{2.7}$ &5(24)&9.19E-03	&6.08E-02	&6.58E-02&	2.43E-01&	3.60&3.60
		\\
		\hline
		$\sqrt{2.8}$ &5(27) &9.14E-03&	6.05E-02&	6.54E-02	&2.43E-01&	3.60&	3.60
		\\
		\hline
		$\sqrt{2.9}$ & 5(27)&9.14E-03	&6.05E-02&	6.54E-02&	2.43E-01&	3.60&	3.60
		\\
		\hline
		$\sqrt{3.0}$ & 5(34)&9.02E-03&	6.01E-02&	6.47E-02&	2.43E-01&	3.60&	3.60
		\\
		\hline
		$\sqrt{3.8}$ & 7(426)&8.79E-03&	5.72E-02	&6.01E-02&	2.43E-01&	3.60&	3.60
		\\
		\hline
	\end{tabular}
\end{table}
\begin{table} [H]\small
	\caption{ Errors of geometric FAS algorithm when $h(H)=\frac{1}{32}(\frac{1}{16})$}\label{32-fas-g}
	\centering
	\begin{tabular}{|c|c|c|c|c|c|c|c|c|  }
		\hline
		$L $ &   iter(iter-c)  &    $\|e_{ u }\|_0$   &  $\|e_{ p }\|_0$ &  $\|e_{  n }\|_0$ & $\|e_{ u }\|_1$   & $\|e_{ p  }\|_1$ & $\|e_{ n }\|_1$  \\
		\hline
		${1}$ &2(7)&2.46E-03&	1.96E-02&	2.09E-02&	1.22E-01&	1.80&	1.80
		\\
		\hline
		$\sqrt{1.5}$ &2(11)&2.05E-03	&2.05E-02&	2.20E-02&	1.22E-01&	1.80&	1.80
		\\
		\hline
		$\sqrt{2.7 }$&  3(42)&2.18E-03&	1.68E-02&	1.84E-02&	1.22E-01&	1.80&	1.80
		\\
		\hline
		$\sqrt{2.8}$ & 3(53)&2.14E-03&	1.70E-02&	1.87E-02&	1.22E-01&	1.80&	1.80
		\\
		\hline
		$\sqrt{3.1}$ & 3(183)&2.04E-03&	1.77E-02&	1.98E-02&	1.22E-01&	1.80&1.80
		\\
		\hline
	\end{tabular}
\end{table}

\vskip0.3cm

Combining the two tables above with Tables \ref{gummel-iter},  \ref{Gummel-16}, \ref{Gummel-32} and \ref{Gummel-TG-2}, we can see
\begin{description}
	\item[(1)] Compared with the Gummel algorithm, the geometric FAS algorithm can calculate larger size of computational domain. For example, the Gummel algorithm does not converge such as when $h=\frac{1}{16}$,  $L=\sqrt{2.8}$, but the geometric FAS algorithm can
	be applied to the cases with larger $L$, such as $L=\sqrt{3.8}$. When the size $L$ is the same, the FAS algorithm can achieve the same accuracy as the Gummel algorithm and converge faster than it, such as the number of iterations of the Gummel algorithm is about 18 times that of the FAS algorithm when $h=\frac{1}{16}$, $L=\sqrt{2.7}$.
\item[(2)] Compared with the TG algorithm, the geometric FAS algorithm has higher accuracy when $L$ is larger. In addition, since TG needs to call linearized  methods such as the Gummel algorithm on the coarse grid level, it can be seen from (1) that the geometric FAS algorithm can calculate a larger computational size than the TG algorithm. For example, the approximate upper limit of the size that the TG algorithm can calculate is about $L=\sqrt{2.7}$ when $h=\frac 1 {32}$, while that is more than $L=\sqrt{3.1}$ for the FAS algorithm.

\end{description}	

Due to the fact that the geometric FAS algorithm usually requires users to provide geometric information of multiple layers of grids and nested grids, it is difficult to apply the geometric FAS algorithm to complex problems. Designing a reasonable algebraic FAS algorithm is an effective way to make up this defect. At present, there have been some preliminary research works in this area, e.g. a coarse grid construction method based on element aggregations is provided in \cite{alloush2022development,birken2019preconditioned} .  Next an effective algebraic FAS algorithm is presented for the PNP equation \eqref{1.1-EAFE}.

\subsubsection{Algebraic FAS algorithm}
In this subsection, we first introduce some notations used in the abgebraic FAS algorithm, then present the algebraic FAS algorithm and apply it to solve PNP equations including a practical problem in ion channel.

For any real matrix $A=(a_{ij})_{n\times n}$, denote by
\begin{itemize}
	\item The dependency set of node $i$
	$$N_i = \{j\neq i~|~a_{ij}\neq 0\}.$$
	
	\item The influence set of node $i$
	$$N_i^T = \{j~|~i\in N_j\}.$$
	
	\item For the given strong and weak threshold $\theta_1$ and row sum parameter $\theta_2$, the strong dependency set $S_i$ of node $i$ is defined as follows:
	\begin{itemize}
		\item[$\diamond$] If $\theta_2 = 1$, or~$\theta_2 \in (0,1)$ and $\bigg{|}\frac{\sum\limits_{j=1}^{n}a_{ij}}{a_{ii}}\bigg{|} \le \theta_2$, then
		\begin{eqnarray*}
			S_i=\{j\in N_i~\big{|}~ a_{ij}\left\{
			\begin{array}{l}
				\geq \theta_1 \max\limits_{k\neq i}\{a_{ik}\},~~a_{ii} < 0 \\
				\leq \theta_1 \min\limits_{k\neq i}\{a_{ik}\},~~a_{ii} \geq 0
			\end{array}
			\right.\}.
		\end{eqnarray*}

		\item[$\diamond$] If $\theta_2 \in (0,1)$ and  ~$\bigg{|}\frac{\sum\limits_{j=1}^{n}a_{ij}}{a_{ii}}\bigg{|} > \theta_2$, then $
		S_i = \emptyset.$

	\end{itemize}
\end{itemize}
\begin{itemize}
	\item Condition C1 of C/F splitting \textsuperscript{\cite{ruge1987algebraic}}:
	\begin{itemize}
		
		\item[$\diamond$] Setting~$i\in F$, $\exists~k\in C$, such that~$k\in S_i$;
		\item[$\diamond$] Setting~$i,j\in F$, and if $j\in S_i$, then~$\exists~ k\in S_i\cap C$, such that~$k\in S_j$.
		
	\end{itemize}
	
	\item   C/F splitting satisfies condition C2: $C$ is a certain maximal independent set of $V$,
	 that is, there is no strong dependency or influence relationship between any two elements in $C$.
\end{itemize}

Next, the algebraic type FAS algorithm is presented by combining the coarsening algorithm, interpolation matrix generation algorithm, and geometric FAS algorithm of  AMG method.

 \begin{algorithm}[H]
\caption{Algebraic type~FAS algorithm}\label{a-FAS}
\begin{description}
  \item[Setup:]~
\begin{description}
  \item[S1.]  Generate single level grid $\mathcal{T}_{h_1}$ and submatrix ~$A^\Phi_1$ of stiffness matrix  ~$A_{1}(\cdot)$ corresponding to the discrete system ~\eqref{fulldis2-EAFE}-\eqref{fulldis-EAFE} on $\mathcal{T}_{h_1}$, $ \bar M_1$ and right hand vector $F_1$.
  \item[S2.]  Get the transfer operator between coarse and fine level.
  \begin{description}
    \item[S2.1] Generate CF array by Algorithm \ref{amg-setup}, denote the number of $C$ points by $n_2^a$, and interpolation matrix ~$P_{1}^{2,a}$ from coase level to fine level.
    \item[S2.2] Get degrees of freedom transfer matrix $R_1^{2,u,a}\in R^{n_{2}^a\times n_1}$, the general element of which is defined by
\begin{equation}\label{Ru2-a}
R_1^{2,u,a}(i,j)=\left\{
         \begin{array}{ll}
           \frac{R_1^{2,a}(i,j)}{ \sum\limits_{j=1}^{n_{1}} R_1^{2,a}(i,j)}, & \hbox{if $R_1^{2,a}(i,j) \ne 0$;} \\
           0, & \hbox{else,}
         \end{array}
       \right.,
\end{equation}
where $ R_1^{2,a}=(P_2^{1,a})^T$.

  \end{description}

\item[S3.] Generate the following submatrix of the coarse grid operator $ A_{2}(\cdot)$
\begin{equation}\label{A2-submatrix}
 A^\Phi_2=R_1^{2,a}A^\Phi_1P_2^{1,a}, ~M_2=R_1^{2,a}M_1P_2^{1,a},
\end{equation}
and give the initial value ~$U_1^{(0)}$, stopping tolerance $\epsilon$ and $\epsilon_c$ between coarse and fine level, and maximum times $M_{iter}$ of iteration .
\end{description}

  \item[Solve:] For any $l\ge0$, the algorithm procedure from ~$U_{1}^{(l)}\rightarrow U_{1}^{(l+1)}$ is as follows:
  \begin{description}

 \item[S1.]{\bf Pre-smoothing}:   Given initial value~$  U_1^{( l)} $, perform $ \nu_{1}\ge 0$ times Gummel iteration to get~$U_1^{ ( fs  )}$.

   \item[S2.]{\bf } Soving the following nonlinear system on the coarse grid by Algorithm \ref{csolve-a}
\begin{equation}\label{ce}
A_2\left(  U_2\right) U_2=\tilde{\tau}_2.
 \end{equation}
Here {
\begin{equation}\label{Atau}
    A_2(\cdot)=\vec R_1^{2,a }A_1 ((\vec R_1^{2,a })^T(\cdot) )(\vec R_1^{2,a })^T,~  \tilde{\tau}_2:=\vec R_1^{2,a }A_1 ((\vec R_1^{2,a })^Ty_2 ) y_2+\vec R_1^{2,a }r_1, 
    \end{equation}}
   where $r_1$ is given by \eqref{r1},
\begin{equation}
\vec R_1^{2,a} =diag(R_1^{2,a}, R_1^{2,a}, R_1^{2,a}),~ y_2 = \vec R_1^{2,u,a}   U_1^{(  f s)},~
\vec R_1^{2,u,a} =diag(R_1^{2,u,a}, R_1^{2,u,a}, R_1^{2,u,a}).
\end{equation}

\item[S3.] {\bf Correction on coarse grid}:
\begin{equation}
  U_1^{ ( bs  )}=  U_{h_{1}}^{( fs)} + ( \vec R_1^{2,a } )^T(  U_2-y_2).
\end{equation}

\item[S4.] {\bf Post-smoothing }: Taking $U_1^{( bs)} $ as the initial value, perform $\nu_2\ge 0$ times Gummel iteration to get $  U_1^{ ( l+1  )} $.
\end{description}
\end{description}
\end{algorithm}

Algorithms ~\ref{amg-setup} and ~\ref{csolve-a} used above are descripted as follows.
\begin{algorithm}[H]\small
\caption{Generate CFmarker array and interpolate matrix $P_{1}^{2,a}$}\label{amg-setup}
\begin{description}
  \item[S1]Given coarsing times $k$, strong and weak threshold $\theta_1$, row sum parameter $\theta_2$, the finest grid matrix $A_1$, generate array~CF\_marker\_l on each level and interpolant matrix $P_l$ between coarse and fine level. For $l=1,\ldots, k$, perform
\begin{description}
  \item[S1.1] Get strong connectivity  matrix  $S_l\in R^{n_l\times n_l}$, the general element of which is
    \begin{equation*}
    S_l(i,j) =\left\{
     \begin{array}{l}
      1,~~j\in S_i^l;\\
      0,~~otherwise.
     \end{array}
    \right.
  \end{equation*}

  \item[S1.2] Let the set of numbers of the degrees of freedoms on level $l$ be $V_l=\{1,2,\cdots,n_l\}$. Perform C/F splitting for $V_l$ by using $S_l$, and get the C/F property array
  \textsuperscript{\cite{ruge1987algebraic}},
   \begin{eqnarray*}
 \text{ CF\_marker\_l}(i) =
  \begin{cases}
    -1, & \mbox{if}~i\in F_l  \\
    1,  & \mbox{if}~i\in C_l
  \end{cases},~C_l\cup F_l=V_l, ~C_l\cap F_l=\emptyset,
\end{eqnarray*}
which strictly satisfies C1 condition, and tries to satisfy~C2 condition \textsuperscript{\cite{ruge1987algebraic}}.

   \item[S1.3]
Get interpolant matrix $P_l\in R^{n_l\times n_{l+1}}$ from the grid~$(l+1)$ to grid $l$, the general element pf which satisfies
\begin{equation*}
  P_l(i,j) = - \frac{1}{a^l_{ii}+\sum\limits_{k\in D_i^{w,l}\cup F^l_i} a^l_{ik}}
  \left( a^l_{ij}
  + \sum\limits_{k\in D_i^{s,l} \backslash F^l_i} \frac{a^l_{ik} \hat{a}^l_{kj}}{\sum\limits_{m\in C^l_i} \hat{a}^l_{km} }\right),
\end{equation*}
where
\begin{equation*}
  \hat{a}^l_{ij} = \left\{
   \begin{array}{ll}
     0, & \text{if sign}(a^l_{ij})= \text{sign} (a^l_{ii}),\\
     a^l_{ij}, & \text{otherwise}.
   \end{array}
  \right.
\end{equation*}

  \item[S1.4] Generate the coefficient matrix on grid~$(l+1)$
  \begin{equation}
  A_{l+1}=P_l^T A_{l}P_l=(a_{ij}^{l+1})_{n_{l+1}\times n_{l+1}}.
  \end{equation}
\end{description}

 \item[S2]  Merge layer by layer for ~CF\_marker\_l to get two-level CFmarker needed and corresponding $n_2^a\times n_1$ prolongation matrix $P_{1}^{2,a}=\Pi_{l=1}^{k} P_l$, where~$n_2^a=n_{k+1}$.
\end{description}
\end{algorithm}

\begin{algorithm}[H]
\small
\caption{Gummel algorithm on coarse system~\eqref{ce}}\label{csolve-a}
\begin{description}
  \item[Setup:] Generate matrix $A_2^{\Phi}$, $\bar M_2 $ and $\tilde{\tau}_2$ from \eqref{A2-submatrix} and~\eqref{Atau}, and give the initial value~$U_2^{(0)}$, stopping tolerance~$\epsilon_c$ and maximum number $M_{iter}$ of iterations.
  \item[Solve:]  For~$l\ge0$,  Gummel iteration from~$U_2^{(l)}\rightarrow U_2^{(l+1)}$ is as follow
\begin{description}
  \item[S1.] Solve the following linear system
  \begin{align}
  A_2^{\Phi}  \Phi_2^{(l+1)}= \tilde{\tau}_2^\Phi-   \bar M_2 P_2^{(l)};
  \end{align}
  \item[S2.]Generate~$A_2^{P} (\Phi_2^{(l+1)})= R_1^{2,a }A_1^{P} ( P_2^{1,a}\Phi_2^{(l+1)})P_2^{1,a }.$
  \item[S3.]Solve the linear system
  \begin{align}
  A_2^{P} (\Phi_2^{(l+1)})P_2^{(l+1)}= \tilde{\tau}_2^{P};
  \end{align}
  \item[S4.] If
\begin{align}\label{tol-3}
\| \tilde{\tau}_2-A_2(  U_2^{ ( l+1 )} )  U_2^{ ( l+1 )}  \|_{R^{3n_{2}^a}}\leq\epsilon_c,
\end{align}
or~$l+1\ge M_{iter}$, then set~$U_1=U_1^{ ( l+1 )}$. Otherwise let~$l\leftarrow l+1$,  return ~S1 and continue.
\end{description}
\end{description}
\end{algorithm}


Next, we compare the geometric FAS algorithm, a coarsening rate which is 1:8, with the algebraic FAS algorithm on two non-uniform grids.

\begin{figure}[H]
	\begin{minipage}[t]{0.5\linewidth}
		\centering
		\includegraphics[height=5.5cm,width=6cm]{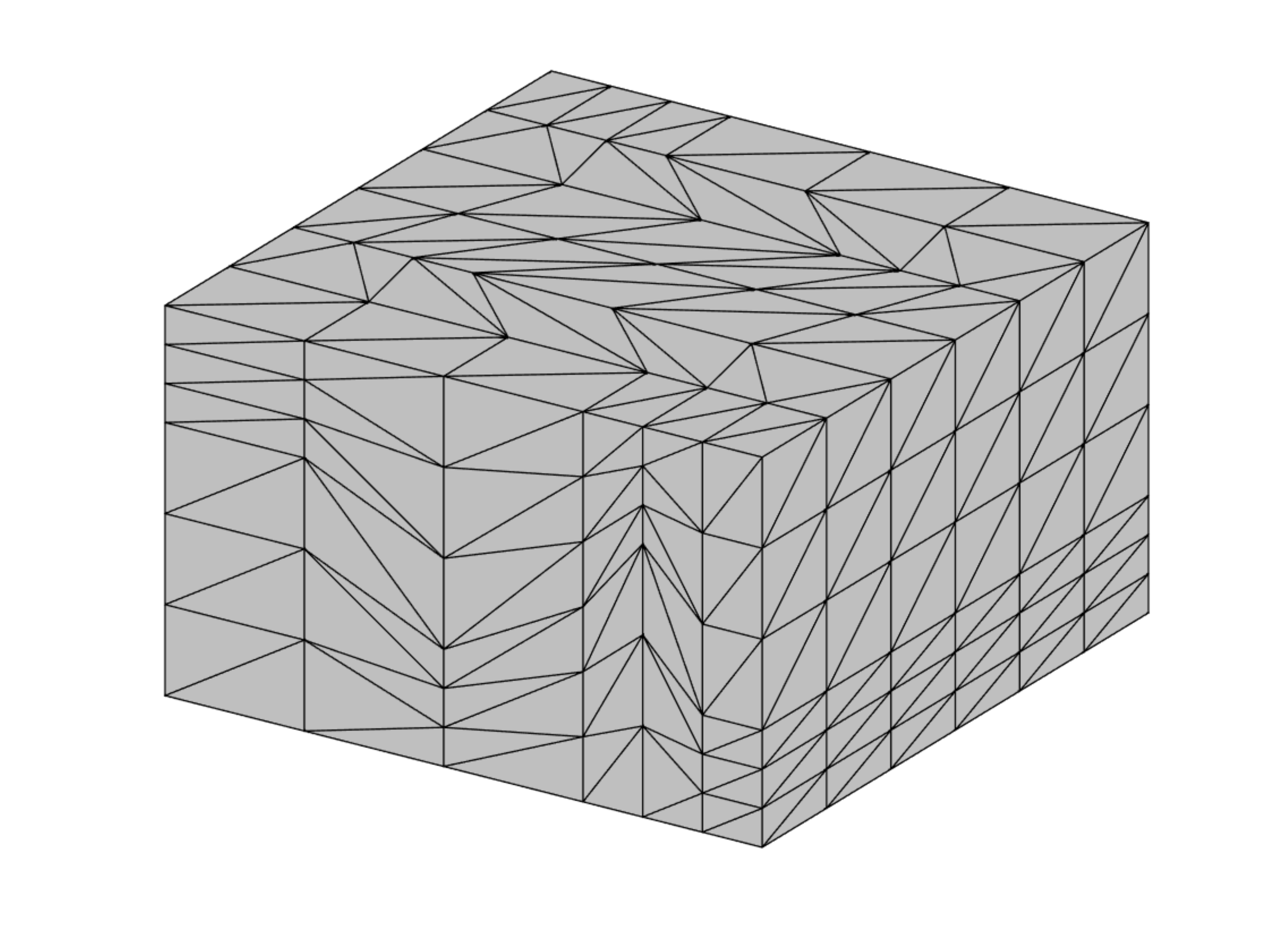}
		\caption{\small $\mathcal {T}_h^1$: Diagram of Kershaw grid with s=0.3 }\label{Kershaw}
	\end{minipage}
	\begin{minipage}[t]{0.65\linewidth}
		\centering
		\includegraphics[height=5.5cm,width=9cm]{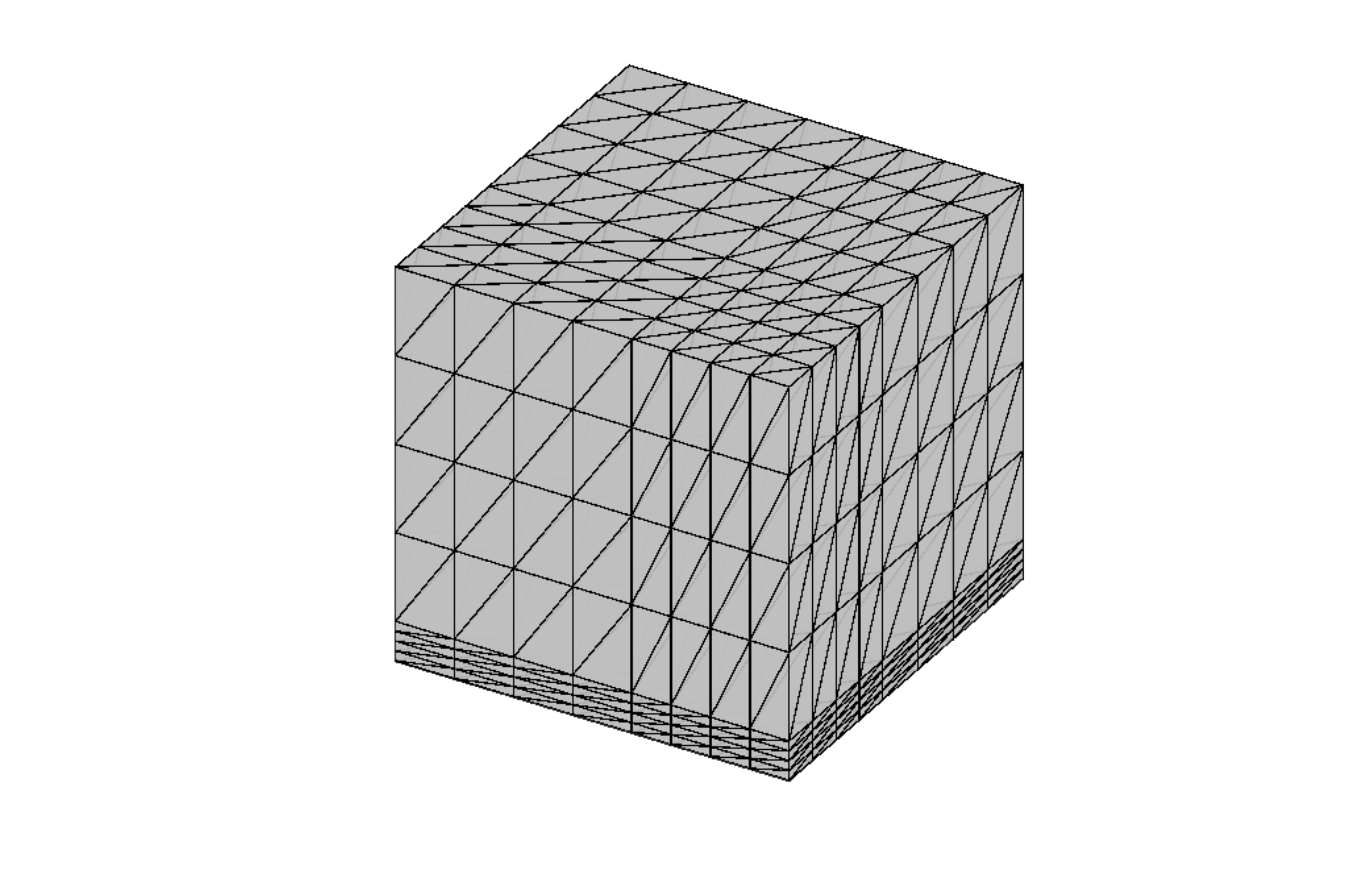}
		\caption{\small $\mathcal {T}_h^2$: Diagram of perturbation Kershaw grid with s=0.2}\label{rand}
	\end{minipage}
\end{figure}

Figure \ref{Kershaw} shows the Kershaw grid \textsuperscript{\cite{kershaw1981differencing}}, where $s=(0,0.5]$ represents the degree of mesh distortion (The worse the quality of mesh becomes as the smaller $s$ is. It is a uniform mesh when $s=0.5$). Figure \ref{rand} shows a perturbation grid, which is obtained  from a uniform grid by perturbing internal nodes as follows
\begin{equation*}
x=\left\{
\begin{array}{ll}
a+(b-a) (1-s) (i-1)/nx, & \hbox{$2 \le i\le  \frac{nx}{2}+1 $ ;} \\
\tilde a_2+(b-\tilde a_2) \frac{ i-\frac{nx}{2}+1 }{nx/2}, & \hbox{$\frac{nx }{2} +1  < i\le  nx $ ;}
\end{array}
\right.
\end{equation*}
\begin{equation*}
y=\left\{
\begin{array}{ll}
a+(b-a) (1-s) (j-1)/ny, & \hbox{$2 \le j\le  \frac{ny}{2}+1 $ ;} \\
\tilde a_2+(b-\tilde a_2) \frac{ j-\frac{ny}{2}+1 }{ny/2}, & \hbox{$\frac{ny }{2} +1  < j\le  ny $ ;}
\end{array}
\right.
\end{equation*}
\begin{equation*}
z=\left\{
\begin{array}{ll}
a+(b-a) s (k-1)/nz, & \hbox{$2 \le k\le  \frac{nz}{2}+1 $ ;} \\
\tilde a_1+(b-\tilde a_1) \frac{ k-\frac{nz}{2}+1 }{nz/2}, & \hbox{$\frac{nz }{2} +1  < k\le  nz $ ;}
\end{array}
\right.
\end{equation*}
where $\tilde a_1=a+\frac{(b-a)s}{2}$, $\tilde a_2=a+\frac{(b-a)(1-s)}{2}$, the left and right endpoints $a=-0.5,~b=0.5$ of the interval, $nx=ny=nz$ are the numbers of partitions in the $x$, $y$, $z$ directions, respectively.

The following two tables display the number of iterations and errors of the corresponding discrete system of the geometric FAS and algebraic FAS algorithms for Example \ref{exam1} on two sets of grids $\mathcal {T}_h^1$ and $\mathcal {T}_h^2$, where $s$ represents the degree of perturbation of grid, iter and iter-c are the iteration times of FAS and the average iteration times of coarse grid, respectively.

\begin{table} [H]
	\caption{ Errors on $\mathcal {T}_h^1$ with  $h=\frac{1}{12}$ and $L=\sqrt{3.5}$ }
	\small\centering
	\begin{tabular}{|c|c|c|c|c|c|c|c|c| }
		\hline
		&s &  iter(iter-c)& $\|e_{ u }\|_0$   &  $\|e_{ p  }\|_0$ &  $\|e_{ n }\|_0$ & $\|e_{ u }\|_1$   & $\|e_{ p  }\|_1$ & $\|e_{ n }\|_1$  \\
		\hline
		geometric FAS  & 0.4  & 13(30)  &2.10E-02&	1.20E-01	&1.66E-01&	3.77E-01&	5.60&  5.57
		\\
		\hline
		algebraic FAS  & 0.4 & 7(41)     &2.10E-02&	1.20E-01	&1.66E-01&	3.77E-01&	5.60&  5.57
		\\
		\hline
		geometric FAS  &0.3 & 12(30)  & 5.10E-02&	3.04E-01&	3.97E-01	&5.98E-01&	8.94& 8.82
		\\
		\hline
		algebraic FAS  & 0.3 & 7(19) & 5.10E-02&	3.04E-01&	3.97E-01	&5.98E-01&	8.94& 8.82
		\\
		\hline
	\end{tabular}
\end{table}

\begin{table} [H]
	\caption{  Errors on $\mathcal {T}_h^2$ with $h=\frac{1}{16}$ and $L=\sqrt{3.0}$  }
	\small\centering
	\begin{tabular}{|c|c|c|c|c|c|c|c|c| }
		\hline
		&s &  iter(iter-c) & $\|e_{ u }\|_0$   &  $\|e_{ p  }\|_0$ &  $\|e_{ n }\|_0$ & $\|e_{ u }\|_1$   & $\|e_{ p  }\|_1$ & $\|e_{ n }\|_1$  \\
		\hline
		geometric FAS  &0.3 & 8(31)  &1.41E-02	&9.58E-02	&1.02E-01	&3.01E-01&	4.46&	4.46
		\\
		\hline
		algebraic FAS  & 0.3 & 7(22)  &1.41E-02	&9.58E-02	&1.02E-01	&3.01E-01&	4.46&	4.46
		\\
		\hline
		geometric FAS  &0.2 &14(29)   &  1.47E-02	&9.87E-02	&1.05E-01	&3.10E-01&	4.59&	4.59
		\\
		\hline
		algebraic FAS &0.2  &6(27)   &  1.47E-02	&9.88E-02	&1.05E-01	&3.10E-01&	4.59&	4.59
		\\
		\hline
	\end{tabular}
\end{table}

\vskip0.3cm

From the above two tables, it can be seen that the geomeric FAS algorithm needs more iterations than algebraic one when $L$ is large, the errors of both under the same $h$ are almost the same, and there are too many iterations on coarse grid for both of them.  For instance, the number of iterations of the geometric FAS algorithm with $s=0.4$ is about twice that of the algebraic one, and numbers of the iteration on the coarse grid of both algorithms are greater than or equal to 30.

Next we present an example of practical PNP equations in ion channel. Since the mesh used in ion channel is unstructured grid, it is difficult to get natrual coarse grid. Hence the geometric FAS algorithm is not suitable and we use algebraic FAS algorithm \ref{a-FAS} to solve the practical PNP problem and compare the result with that of Gummel iteration from the error, convergence speed and CPU time.

\begin{figure}[H]\centering
\includegraphics[width=5cm,height=6cm]{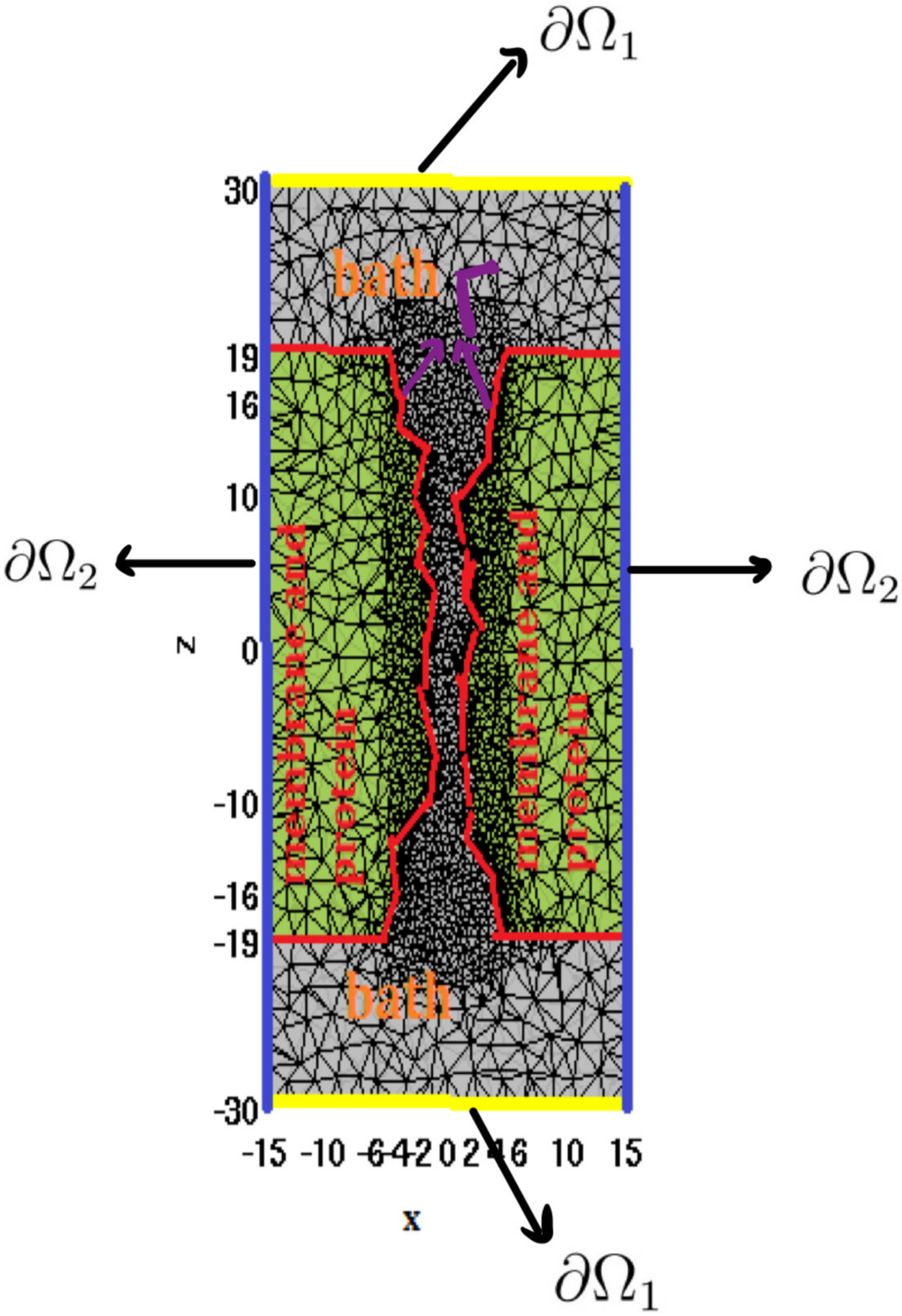}
\caption {2D cut plane y=0 for the mesh of gramicidin A single channel, where the green and grey parts shows the solute $\Omega_m$ and solvent domain $\Omega_s$ resepectively} \label{cutplanebiid}
\end{figure}

\begin{exam}
Consider the following PNP equations for simulating Gramicidin A(gA) ion channel in  1:1  CsCl solution with valence $+1$ and $-1$
\textsuperscript{\cite{lu2007electrodiffusion}}:
\begin{equation}\label{prob1}
\left\{\begin{array}{ll}
-\nabla\cdot\Big(-  D_i\left(\nabla
p^i+\frac{e}{K_BT}q^ip^i\nabla\phi\right) \Big)=0, \text{in}~\Omega_s,~i=1,2,\\
-\nabla\cdot(\epsilon\nabla\phi)=\lambda\sum\limits_{i=1}^2q^ip^i +\rho^f,~\text{ in}~\Omega=\Omega_s\cup\bar\Omega_m,
\end{array}\right.
\end{equation}%
where the unknown function $\phi$ is electrostatic potential and $p^1$ and $p^2$ are the concentrations of positive and negative ion, respectively. Here  $\Omega=[-15, 15]\times[-15,15]\times [-30,30] {\AA}$, $\Omega_s$ represents the bulk regions,   $\Omega_m$ represents the membrane and protein region, the diffusion coefficient ~$D_{1} =2.0561 \times 10^{-9} \mathrm{~m}^{2}/\mathrm{s}$, $D_{2}=2.0321\times 10^{-9} \mathrm{~m}^{2} / \mathrm{s} $,  charges $q^1=1$ and $q^2=-1$, elementary charge $e =1.6 \times 10^{-19} \mathrm{C}$, Boltzmann constant $ K_{B}=4.14 \times 10^{-21} \mathrm{J}$, the dielectric coefficient $\epsilon(x)=\left\{\begin{array}{ll} 2\epsilon_0\ ,~ \mbox{in}~ \Omega_m\\
80\epsilon_0 ,~ \mbox{in}~ \Omega_s
\end{array}\right.$, with permittivity of vacuum~$\varepsilon_{0}=8.85 \times 10^{-12} \mathrm{C}^{2}/\left(\mathrm{N} \mathrm{m}^{2}\right)$,  $\lambda=\left\{\begin{array}{ll} 0,~ \mbox{in}~\Omega_m\\
1,~ \mbox{in}~\Omega_s
\end{array}\right.$, $\rho^f({\bf x})=\sum\limits_j^{N_m} q_j\delta(  {\bf x}-  {\bf x_j})$ is an ensemble of singular atomic charges $q_j$ located at $x_j$ inside the protein.


Let $\Gamma=\partial\Omega_s\cap\partial\Omega_m$ be the interface that intersect the membrane and the bulk region, see the red curve in Figure~\ref{cutplanebiid}, $\partial\Omega_1$ is the part of boundary of $\Omega$ that perpendicular to $z$ axis (suppose the channel is along the $z$ axis), $\partial\Omega_2$ is the part of boundary of $\Omega$ that is along the channel, 
$\partial\Omega_3 =\partial\Omega \setminus \Gamma $ is the surface of the bulk region from which the interface of membrane $\Gamma$ is removed.

The interface and boundary conditions are
\begin{equation}
[\phi]=0,~[\epsilon\frac{\partial\phi}{\partial \upsilon^m}]=0,~~\mbox{on}~~\Gamma \end{equation}
and
\begin{equation}\label{1bd}\left\{\begin{array}{ll}
\phi=V_{applied},~\mbox{on}~\partial\Omega_1,\\
\frac {\partial\phi}{\partial\upsilon }=0,~\mbox{on}~ \partial\Omega_2 ,\\
D_i\left(\nabla
p^i+\frac{e}{K_BT}q^ip^i\nabla\phi\right) \cdot \upsilon =0,~\mbox{on}~\Gamma,~i=1,2,\\
p^i=p_{\infty},~\mbox{on}~\partial\Omega_3,~i=1,2,
\end{array}\right.\end{equation}
where $V_{applied}$ is the applied potential,  $\upsilon^m$ and  $\upsilon$ are the unit outer normal vectors of $\Omega_m$ and the boundary, respectively,  and $p_{\infty}$  is the initial concentration.

\end{exam}
In the simulation,  $\mathcal {T}_h$ is a tetrahedral mesh generated using TMSmesh \textsuperscript{\cite{chen2011tmsmesh}}, including 224650 tetrahedral elements and 37343 nodes. The number $N_m$ of atoms , charge number $q_j$, and the position $x_j$ of the $j$th atom in proteins are obtained from the protein data bank.

Since it is difficult to find the analytical solution of \eqref{prob1}-(\ref{1bd}), the result of simulation is compared with the experiment result via the current to observe  the accuracy of the algorithms. The electric current is as follows \textsuperscript{\cite{tu2013parallel,lu2011poisson}}
\begin{equation}
I=e\int_{\Omega_I}(J_1-J_2)dS,
\end{equation}
where ${\Omega_I}$ is any cross section of the channel, and the fluxes of positive and negative ions are defined by
\begin{equation}
J_i=-\nabla \cdot D_i\left(\nabla
p^i+\frac{e}{K_BT}q^ip^i\nabla\phi\right) ,i=1,2.
\end{equation}

The following two tables show the currents, and number of iterations and CPU time respectively by using Gummel iteration and algebraic FAS algorithm under different voltages. The experimental data was provided by Andersen\textsuperscript{ \cite{andersen1983ion}}.

\begin{table}[H]\centering\small
\caption{The current from the numerical simulation and experimental data at CsCl
	concentration 0.02M and different voltage(mV)
	}
\label{ion-A}
\vskip0.3cm
\begin{tabular}{|c|c|c|c|c|c|}
\hline
{\small Voltage (mV)} &  {\small  Experimental
	data(pA) } &  {\small  Gummel results (pA) }  &  {\small FAS results (pA)}    \\
\hline
25 & 0.13 & 0.14  &    0.14 \\
50 & 0.24 & 0.26   &   0.26\\
75 & 0.34 & 0.34 &    0.34\\
100 & 0.40 & 0.40  &   0.40\\
150 & 0.47 & 0.50   &  0.50\\
200 & 0.53 & 0.57   &  0.57 \\
250 & 0.56 &  0.65  &  0.65\\
300 & 0.60 &  0.73 &    0.73\\
350 & 0.62 &  0.81  &   0.81\\
400 & 0.65 &  0.89 &    0.89\\
\hline
\end{tabular}
\end{table}
\begin{table}[H]\centering\small
\caption{ Number of external iterations and the CPU time(s) for Gummel and FAS algorithms at CsCl concentration
	0.02M and different voltage(mV)
	}
\label{ion-B}
\vskip0.3cm
\begin{tabular}{|c|c|c|c|c|c|}
\hline
{\small Voltage(mV) }& {\small   Gummel iteration number} &  {\small FAS iteration number }&  {\small Gummel CPU time }&{\small  FAS CPU time } \\
\hline
25 & 139 & 1 & 690 & 29  \\
50 &  139 & 1 & 688 & 29  \\
75 &  139  & 1 & 688 & 29 \\
100 & 139  & 1 & 679 & 29 \\
150 & 145  & 1 & 717 & 28 \\
200 & 145 & 1 & 722 & 28 \\
250 & 145 & 1 & 717 & 28 \\
300 & 145 & 1 & 725 & 28 \\
350 & 145 & 1 & 721 & 29 \\
400 & 145 & 1 & 718 & 28 \\
\hline
\end{tabular}
\end{table}
\begin{remark}
One FAS iteration includes two rounds of pre-smoothing and post-smoothing respectively, and solving a coarse grid system once. Therefore, the time required to perform a FAS iteration is approximately six times that of a Gummel iteration.
\end{remark}

From Tables \ref{ion-A} and \ref{ion-B}, we can see FAS iteration runs much faster than Gummel iteration and the accuracy of them are similar.

Although the algebraic FAS method has better adaptability on unstructured or non-uniform grids, there are still defects such as insufficient adaptability to $L$ and excessive iteration for solving corresponding coarse grid nonlinear algebraic systems. Next, further improvements will be made based on the acceleration technique and adaptive method.
\vspace{3mm}

\subsection{Improved algorithms based on the acceleration techniques and adaptive methods}

First, we discuss some improved techniques for Gummel iteration. A natural idea is to use the under-relaxation technique \textsuperscript{\cite{nonner2001ion,lu2007electrodiffusion}},
\begin{eqnarray}
u^{new}=\alpha u^{new }+(1-\alpha)u^{old},
\end{eqnarray}
where the relaxation parameter $\alpha\in(0,1)$ is a preset constant. The corresponding Gummel iteration is described as follows.
\begin{algorithm}[!h]
	\caption{Gummel iteration based on under-relaxation technique}\label{gummel-dichotomy}
	\begin{description}
		\item[S1.] Initialization.  Given initial values $\phi_h^{(0)}=0,~p_h^{1,(0)}=0,~p_h^{2,(0)}=0$, iteration time $k=0$, parameter $\alpha$, tolerance $tol=1E-6$, maximum iteration number $M_{iter}=1000$.
		
		\item[S2.] Solve the linear equation
		\begin{equation}\label{disconh}
		(  \nabla \phi_h^{(k+1)} ,  \nabla v_h)=(f+ p_h^{1,(k)} - p_h^{2,(k)} ,v_h) ,~  \forall v_h \in V_h
		\end{equation}
	to get $\phi_h^{(k+1)}$, and then let $\phi_h^{(k+1)}=\alpha \phi_h^{(k+1)}+(1-\alpha) \phi_h^{(k)}$.
		
		\item[S3.] Solve the linear equations
		\begin{equation}\label{disconH}
		\left\{
		\begin{array}{ll}
		(  \nabla   p_h^{1,(k+1)} ,\nabla v_h)+  (     p_h^{1,(k+1)} \nabla  \phi_h^{(k+1)} ,\nabla v_h)= (F^1,v_h),
		&\forall v_h\in V_h  \\
		(  \nabla   p_h^{2,(k+1)}  ,\nabla v_h)-  (     p_h^{2,(k+1)}  \nabla  \phi_h^{(k+1)}  ,\nabla v_h)=  (F^2,v_h),
		&\forall v_h \in V_h
		\end{array}
		\right.\end{equation}
		to get $p_h^{1,(k+1)}$ and $p_h^{2,(k+1)}$, and set $ p_h^{1,(k+1)}=\alpha  p_h^{1,(k+1)}+(1-\alpha ) p_h^{1,(k)}$, $ p_h^{2,(k+1)}=\alpha p_h^{2,(k+1)}+(1-\alpha )p_h^{2,(k)}$.
		
		\item[S4.] If $\|\phi_h^{(k+1)} -\phi_h^{(k)} \|_{L^2}\ge tol $ and $k\le M_{iter}$, then $k=k+1$, and return to {\bf S2}; otherwise output the numerical solution
		$$\phi_h=\phi_h^{(k+1)},~p_h^{1}=p_h^{1,(k+1)},~p_h^2=p_h^{2,(k+1)}.$$
	\end{description}
\end{algorithm}

Next, we consider Example \ref{exam1} to investigate the adaptability to convective intensity and convergence speed of the above algorithm. Set $\alpha=0.5$, $tol=1E-6$. Table \ref{gummel-dichotomy-iter} shows the numbers of iteration on uniform grids.
\begin{table} [H]
	\caption{Numbers of iterations}\label{gummel-dichotomy-iter}
	\centering
	\begin{tabular}{|c|c|c|c|c|c|c|c|c|c|c|c|c| }
		\hline
		\diagbox{$h$}{$L$}  &  $\sqrt{2.6 }$ &  $\sqrt{2.7}$    &  $\sqrt{14.0}$ &  $\sqrt{15.0}$   \\
		\hline
		$\frac{1}{16}$ &20 & 20  &21 &*
		\\
		\hline
		$\frac{1}{32}$ &19  & 19   & 19 &*
		\\ \hline
	\end{tabular}
\end{table}
where * represents divergence.

Comparing Table \ref{gummel-dichotomy-iter} with Table \ref{gummel-iter}, it can be seen that Algorithm \ref{gummel-dichotomy} indeed enhances its adaptability to $L$ (i.e. convective intensity), and the number of iterations does not change significantly with the scale, but there may be sudden instances of divergence. It may due to the fixation of values of $\alpha$. Hence, the following accelerated Gummel algorithm is provided to solve this problem.

\subsubsection{Acceleration algorithms}
Consider the following nonlinear algebraic system:
\begin{align}\label{A(U)UF-g}
A(U)U=F,
\end{align}
where $A(U)$ and $F$  are order $n$ square matrix and vector, respectively.

Let $\hat U^k$ be the present iteration vector by using Gummel iteration to solve \eqref{A(U)UF-g} and define the corresponding minimum problem
\begin{align}\label{minr}
\min_{\alpha\in [0,1] }\|\alpha\hat r^k+(1-\alpha) r^{k-1}\|^2_{R^n},
\end{align}
where
\begin{equation}\label{rkrk-1}
\hat r^k=F-A(\hat U^k)\hat U^k, ~~r^{k-1}=F-A(U^{k-1})U^{k-1}.
\end{equation}
It is easy to know the minimum point is
\begin{equation}\label{alpha*}
\alpha^*=-\frac{b}{2a},~
a= ((\hat r^k)^T , \hat r^k )+ (     (r^{k-1})^T,      r^{k-1})- 2( (\hat r^k)^T , r^{k-1}) \ge0,
\end{equation}
\begin{equation}\label{alpha*-b}
~b=2 \big(( (\hat r^k)^T ,  r^{k-1})- (   (r^{k-1})^T,  r^{k-1})\big).
\end{equation}

Updating $\hat U^k$ with $U^k$ using the parameter $\alpha^*$, we have
\begin{eqnarray}\label{u-update}
U^k=\hat U^k+(1-\alpha^*)\hat U^{k-1}.
\end{eqnarray}
From the minimum problem \eqref{minr} and \eqref{u-update}, we get the following accelerated Gummel iteration.

\begin{algorithm}[!h]
	{
		\caption{Accelerated Gummel iteration-I}\label{Gummel-js-0}
		\begin{description}
			\item[S1.] Initialization. Given initial value $\phi_h^{(0)}=0,~p_h^{1,(0)}=0,~p_h^{2,(0)}=0$, number of iterations $k=0$, stopping tolerance $tol=1E-6$, maximum iteration number $M_{iter}$=1000.
			\item[S2.] Solve linear equations
			\begin{equation}\label{disconh}
			(  \nabla \phi_h^{(k+1)} ,  \nabla v_h)=( f+p_h^{1,(k)}-p_h^{2,(k)} ,v_h) ,~~ \forall v_h \in V_h
			\end{equation}
			to get $\phi_h^{(k+1)}$.
			
			\item[S3.] Solve the following linear equations
			\begin{equation}\label{disconH}
			\left\{
			\begin{array}{ll}
			(  \nabla   p_h^{1,(k+1)} ,\nabla v_h)+  (     p_h^{1,(k+1)} \nabla  \phi_h^{(k+1)} ,\nabla v_h)= (F^1,v_h),
			&\forall v_h\in V_h  \\
			(  \nabla   p_h^{2,(k+1)}  ,\nabla v_h)-  (    p_h^{2,(k+1)}  \nabla  \phi_h^{(k+1)}  ,\nabla v_h)=  (F^2,v_h),
			&\forall v_h \in V_h
			\end{array}
			\right.\end{equation}
			to obtain $p_h^{1,(k+1)}$ and $p_h^{2,(k+1)}$.
			
			\item[S4.] If $k\ge 2$, 
			then compute $$r_\phi^{(k+1)}= (  f+p_h^{1,(k+1)} - p_h^{2,(k+1)} , \Psi) -(  \nabla \phi_h^{(k+1)} ,  \nabla\Psi). $$
			From the minimum problem \eqref{minr}, we get $\alpha^* $,  and then let $\phi_h^{(k+1)}=\alpha ^* \phi_h^{(k+1)}+(1-\alpha^* )\phi_h^{(k)}$, $p_h^{1,(k+1)}=\alpha ^* p_h^{1,(k+1)}+(1-\alpha^* )p_h^{1,(k)}$, $p_h^{2,(k+1)}=\alpha ^* p_h^{2,(k+1)}+(1-\alpha^* )p_h^{2,(k)}$.
			
			\item[S5.] If $|\phi_h^{(k+1)} -\phi_h^{(k)} \|_{L^2}\ge tol $ and $k\le M_{iter}$ then set $k=k+1$ and return to S2, otherwise output the numerical solution
			$$\phi_h=\phi_h^{(k+1)},~p_h^{1}=p_h^{1,(k+1)},~p_h^{2}=p_h^{2,(k+1)}.$$
		\end{description}
}\end{algorithm}
{
	\vskip0.3cm

Next, we investigate the adaptability to convective intensity and convergence speed of Algorithm \ref{Gummel-js-0} by taking Example \ref{exam1} as an example. Let $tol=1E-6$. The following Table \ref{gummel-iter-js-0} shows the numbers of iterations on uniform meshes.
\begin{table} [H]
		\caption{Numbers of iteration of Algorithm \ref{Gummel-js-0}}\label{gummel-iter-js-0}
		\centering{
			\begin{tabular}{|c|c|c|c|c|c|c|c|c|c|c|c|c| }
				\hline
				\diagbox{$h$}{$L$}  &  $\sqrt{2.6 }$ &  $\sqrt{2.7}$    &  $\sqrt{7.4}$ &  $\sqrt{7.5}$&  $\sqrt{7.6}$&  $\sqrt{7.7}$     \\
				\hline
				$\frac{1}{16}$ & 9	& 9   &  16 &  18 & 21 & --
				\\
				\hline
				$\frac{1}{32}$ & 9	& 9  &18 & -- & -- & --
				\\ \hline
		\end{tabular}}
	\end{table}
	Comparing Table \ref{gummel-iter-js-0} with Table \ref{gummel-dichotomy-iter}, we can see that the number of itertion of Algorithm \ref{Gummel-js-0} is less than than that of Algorithm \ref{gummel-dichotomy} in the convergent cases, but the adaptability to $L$ of it is worse than that of Algorithm \ref{Gummel-js-0}.

Note that $\alpha^*$ is obtained by solving the minimum problem with respect to the residue of the Poisson equation. Hence, it may not be a good parameter for $p_h^1$ and $p_h^2$. To solve this problem, the combined $\phi_h^{(k+1)}$ is used to update $p_h^{1,(k+1)}$ and $p_h^{2,(k+1)}$. The detailed algorithm is as follows.	

}

\begin{algorithm}[!h]
	\caption{ Accelerated Gummel iteration-II}\label{Gummel-js}
	\begin{description}
		\item[S1.] Initialization. Given initial value $\phi_h^{(0)}=0,~p_h^{1,(0)}=0,~p_h^{2,(0)}=0$, number of iteration $k=0$, stopping tolerace $tol=1E-6$, maximum itertion number $M_{iter}$=1000.
		
		\item[S2.] Solve the linear equation
		\begin{equation}\label{disconh}
		(  \nabla \phi_h^{(k+1)} ,  \nabla v_h)=( f+p_h^{1,(k)} - p_h^{2,(k)} ,v_h) ,~~ \forall v_h \in V_h
		\end{equation}
		to get $\phi_h^{(k+1)}$.
		
		\item[S3.] Solve the linear equations
		\begin{equation}\label{disconH}
		\left\{
		\begin{array}{ll}
		(  \nabla   p_h^{1,(k+1)} ,\nabla v_h)+  ( p_h^{1,(k+1)} \nabla  \phi_h^{(k+1)} ,\nabla v_h)= (F^1,v_h),
		&\forall v_h\in V_h,  \\
		(  \nabla   p_h^{2,(k+1)}  ,\nabla v_h)-  (p_h^{2,(k+1)}  \nabla  \phi_h^{(k+1)}  ,\nabla v_h)=  (F^2,v_h),
		&\forall v_h \in V_h
		\end{array}
		\right.\end{equation}
		to obtain $p_h^{1,(k+1)}$ and $p_h^{2,(k+1)}$.
		
		\item[S4.] If $k\ge 2$, then compute
		$$r_\phi^{(k+1)}= (f+ p_h^{1,(k+1)} - p_h^{2,(k+1)} , \Psi) -(  \nabla \phi_h^{(k+1)} ,  \nabla\Psi), $$
	and solve the minimum problem \eqref{minr} to get $\alpha^* $. Set
		$\phi_h^{(k+1)}=\alpha ^* \phi_h^{(k+1)}+(1-\alpha^*)\phi_h^{(k)}$,  and repeat Step {\bf S3} to obtain $n_h^{(k+1)},p_h^{(k+1)}$.
		
		\item[S5.] If $\|\phi_h^{(k+1)} -\phi_h^{(k)} \|_{L^2}\ge tol $ and $k\le M_{iter}$, then set $k=k+1$, return to S2, otherwise output the following solution
		$$\phi_h=\phi_h^{(k+1)},~p_h^{1}=p_h^{1,(k+1)},~p_h^{2}=p_h^{2,(k+1)}.$$
	\end{description}
\end{algorithm}

Now, we still use Example \ref{exam1} to investigate the adaptability to convective intensity and convergence speed of Algorithm \ref{Gummel-js}. Table \ref{gummel-iter-js} shows the numbers of iterations on uniform meshes.

\begin{table} [H]
	\caption{Numbers of iteration of Algorithm \ref{Gummel-js}}\label{gummel-iter-js}
	\centering
	\begin{tabular}{|c|c|c|c|c|c|c|c|c|c|c|c|c| }
		\hline
		\diagbox{$h$}{$L$}  &  $\sqrt{2.6 }$ &  $\sqrt{2.7}$   &  $\sqrt{10}$ &  $\sqrt{11}$&  $\sqrt{12}$&  $\sqrt{13}$&  $\sqrt{14}$ \\
		\hline
		$\frac{1}{16}$ &3	&3   &19&39 &110 & 398  &$>$1000
		\\
		\hline
		$\frac{1}{32}$ &3	&3   &19&43&132&523&$>$1000
		\\ \hline
	\end{tabular}
\end{table}

Comparing Table \ref{gummel-iter-js} with Table \ref{gummel-iter}, we know that Algorithm \ref{Gummel-js} indeed accelerates convergence speed and can calculate larger $L$. For example, the approximate upper limit of $L$ that the Gummel iteration can calculate is about $L=\sqrt{2.7}$ when $h=\frac{1}{16}$, while that is more than $L=\sqrt{13}$ for Algorithm \ref{Gummel-js}, and the number of iteration of Gummel algorithm is about $39$ times that of Algorithm \ref{Gummel-js} when $L=\sqrt{2.6}$.

The above acceleration techniques have improved the adaptability to size $L$ in a certain extent, but when $L$ increases to a certain value, there are also a rapid increase in iteration numbers and slow convergence speed. It is observed that $\alpha^*$ is very small at the first iteration and becomes larger as the decline of the residual. For instance, in the case of $L=\sqrt{11}$, the magnitude of $\alpha^*$ is $1E-3$ when the number of iteration is smaller, while it is $1E-1$ at the last few steps. Considering this phenomenon, we present a adaptive Gummel iteration as follows.

\subsubsection{An adaptive algorithm}
\begin{algorithm}[H]
	\caption{Adaptive Gummel iteration}
	\label{Gummel-js-a}
	\begin{description}
		\item[S1.] Setting $M_{iter}=2$, and given $\theta_\alpha$ and $\alpha$ , call Algorithm \ref{Gummel-js} to get $\alpha^*$.
		
		\item[S2.] If $\alpha^*\ge \theta_\alpha$, then use Algorithm \ref{Gummel-js} directly to solve the algebraic system, otherwise go to S3.
		
		\item[S3.-S4.] Same as Steps S1-S2 in Algorithm \ref{Gummel-js}.
		
		\item[S5.] If $k=1$,  then set $\phi_h^{(k+1)}=\alpha \phi_h^{(k+1)}+(1-\alpha) \phi_h^{(k)}$, and solve linear equation
		~\eqref{disconH}
		 to get $p_h^{1,(k+1)}$ and $p_h^{2,(k+1)}$, otherwise compute
		\begin{eqnarray}
		&r_\phi^{(k+1)}= ( p_h^{1,(k+1)}-p_h^{2,(k+1)} ,\nabla\Psi) -(  \nabla \phi_h^{(k+1)},\nabla\Psi),
		\\& uc=
		\frac{\| \phi_h^{(k+1)}-\phi_h^{(k)}\|_{L^2}}{\| \phi_h^{(k)}\|_{L^2}},
		\end{eqnarray}
		and obtain $\alpha^*$ from the minimum problem \eqref{minr}. Then given $(\theta_1^i, \theta_2^i),i=1,\dots,m$, call Algorithm \ref{updata-phi} to update $\phi_h^{(k+1)}$,  $p_h^{1,(k+1)}$ and $p_h^{2,(k+1)}$.
				
		\item[S6.] Compute
		$$r^{(k+1)}=(r_\phi^{(k+1)},r_{p^1}^{(k+1)},r_{p^2}^{(k+1)}),
		$$
		where    $$ r_{p_h^{i}}^{(k+1)}=(F^i,v_h)-(  \nabla    p _h^{i,(k+1)}  ,\nabla v_h)+ (     p_h^{i,(k+1)}  \nabla  \phi_h^{(k+1)}  ,\nabla v_h),~i= {1}, {2}.$$
		If $\|r^{(k+1)}\|_{L^2}\ge tol $ and $k\le M_{iter}$, then $k=k+1$ and return to S3, otherwise	$$\phi_h=\phi_h^{(k+1)},~p_h^{2}=p_h^{2,(k+1)},~p_h^{1}=p_h^{1,(k+1)}.$$
	\end{description}
\end{algorithm}

\begin{algorithm}[H]
	\caption{Algorithm of updating $\phi_h^{(k+1)}, ~p_h^{1,(k+1)},~p_h^{2,(k+1)}$}
	\label{updata-phi}
	
	\begin{description}
		\item [S1. Update $\phi_h^{(k+1)}$].
		
		\begin{description}
			\item[S1.1]   If $\|r^{(k)}\|_{L^2}\in(r_1, +\infty)$ and $uc\in [\theta_1^1 ,\theta_2^1]$, goto S2, otherwise
			\begin{description}
				\item[S1.1] Setting $\hat \phi_h^{(k+1)}=\alpha^* \phi_h^{(k+1)}+(1-\alpha^*)\phi_h^{(k)}$,
				compute $\hat {uc}=
				\frac{\|\hat \phi_h^{(k+1)}-\phi_h^{(k)}\|_{L^2}}{\| \phi_h^{(k)}\|_{L^2}}$.
				\item[S1.2]~
				\begin{description}
					\item If $\hat {uc}<1e-10$, then take $\alpha^*= 0.5$ and return to S1.1.
					
					\item If $1e-10\le \hat {uc}<\theta_1^1$, then $\alpha^*= 2 {\alpha}^* $ and return to S1.1.
					
					\item If $\hat {uc}\in [\theta_1^1  ,\theta_2^1]$, then set $ \phi_h^{(k+1)}=\hat \phi_h^{(k+1)}$, goto S2.
					
					\item If $\hat {uc}> \theta_2^1$, then set $\alpha^*= {\alpha^*}/{2}$, and return to S1.1.
				\end{description}
			\end{description}
			\qquad $\cdot$
			
			\qquad $\cdot$
			
			\qquad $\cdot$
			\item[S1.m] If $\|r^{(k)}\|_{L^2}\in (r_{m},r_{m-1}]$  and $uc\in [\theta_1^{m } ,\theta_2^{m }]$, then goto S2, otherwise repeat steps S1.1-S1.2, in which  $\theta_1^1$ and $ \theta_2^1$ are considered as $\theta_1^{m}$ and $\theta_2^{m}$, respectively.
		\end{description}
		
		\item[S2. (Update $p_h^{1,(k+1)},p_h^{2,(k+1)}$)]
		 Solve \eqref{disconH} to get $p_h^{1,(k+1)}$ and $p_h^{2,(k+1)}$.
		
	\end{description}
	
\end{algorithm}

Next we use Example \ref{exam1} as an example to observe the adaptability to convective intensity and convergence speed of  algorithm \ref{Gummel-js-a}. The parameters in Algorithm \ref{Gummel-js-a} are taken as follows:
{\small
	\begin{equation}\label{zsy-p-1}
	\begin{aligned}
	&tol=1E-6, ~\alpha=0.1,~\theta_\alpha=1e-3, ~m=3,~ r_1=1e-3,~ r_2=1e-4,~r_3=1e-6,
	\\
	&(\theta_1^1, \theta_2^1)=(1e-1,5e-1), ~(\theta_1^2, \theta_2^2)=(1e-2,1e-1),~(\theta_1^3, \theta_2^3)=(0,100).
	\end{aligned}
	\end{equation}}
Table \ref{gummel-a-iter-1} shows the numbers of iteration of Algorithm \ref{Gummel-js-a} on uniform meshes.
\begin{table} [H]
	\caption{numbers of iteration of Algorithm \ref{Gummel-js-a}}\label{gummel-a-iter-1}
	\centering \small
	\begin{tabular}{|c|c|c|c|c|c|c|c|c|c|c|c|c|c|c| }
		\hline
		\diagbox{$h$}{$L$}  &  $\sqrt{2.7}$  &  $\sqrt{2.8}$  &  $\sqrt{10}$ &  $\sqrt{11}$&  $\sqrt{12}$&  $\sqrt{13}$&  $\sqrt{14}$ &  $\sqrt{15}$ &  $\sqrt{18}$ &  $\sqrt{19}$ &  $\sqrt{20}$&  $\sqrt{25}$  \\
		\hline
		$\frac{1}{16}$ 	& 3 & 3   &10 & 12 & 12 & 13 & 14 &  16  &  17 & -- & -- & --
		\\
		\hline
		$\frac{1}{32}$	&3 &3  & 17  & 16  & 13  & 17 & 24 & 30 &  29 & 27 & 26& --
		\\ \hline
	\end{tabular}
\end{table}
It can be seen by comparing the results of Tables \ref{gummel-iter-js}-\ref{gummel-a-iter-1} that the
adaptability to $L$ of Algorithm \ref{Gummel-js-a} is better than that of Algorithm \ref{Gummel-js}. For example, the approximate upper limit of $L$ that Algorithm \ref{Gummel-js} can calculate is about $L=\sqrt{14}$ when $h=\frac{1}{16}$, while that is more than $L=\sqrt{18}$ for Algorithm \ref{Gummel-js-a}, and the number of iterations of it is not more than that of Algorithm \ref{Gummel-js}. In order to further improve the adaptability to $L$, we adjust the values of parameters as follows:
{\small
	\begin{equation}\label{zsy-p-2}
	\begin{aligned}
	&tol=1E-6, ~\alpha=0.1,~\theta_\alpha=1e-3,
	\\&~m=3,~ r_1=1e-2, ~  r_2=1e-3,~r_3=1e-5,~r_4=1e-6,
	\\&(\theta_1^1, \theta_2^1)=(1e-1,5e-1),~(\theta_1^2, \theta_2^2)=(1e-2,1e-1),~(\theta_1^3, \theta_2^3) =(1e-4,1e-3)
	\\&(\theta_1^4, \theta_2^4)=(1e-5,1e-4),
	\end{aligned}
	\end{equation}}
and show the results on Table \ref{gummel-a-iter-2}.
\begin{table} [H]
	\caption{Numbers of iterations of Algorithm \ref{Gummel-js-a} after parameter adjustment} \label{gummel-a-iter-2}
	\centering
	\begin{tabular}{|c|c|c|c|c|c|c|c|c|c|c|c|c|c|c| }
		\hline
		\diagbox{$h$}{$L$}  &    $\sqrt{18}$ &  $\sqrt{19}$  & $\sqrt{20}$&  $\sqrt{25}$ &  $\sqrt{30}$ &  $\sqrt{35}$ &  $\sqrt{40}$  \\
		\hline
		$\frac{1}{16}$ &116 & 113 & 108&  47  &  151 &  153 &  --
		\\\hline
	\end{tabular}
\end{table}

It is seen from Tables \ref{gummel-a-iter-1} and  \ref{gummel-a-iter-2} that the size $L$ that Algorithm \ref{Gummel-js-a} can calculate can be larger if $m$ is increased and $ \theta_1^\eta, \theta_2^\eta$ is decreased.

It should be pointed out that the improved algorithms based on the acceleration technique and the adaptive method
in this section are designed for the Gummel iterations an example. However, these algorithms essentially do not require modifications and can be directly applied to the FAS algorithm.

\section{Conclusion}
In this paper, we study several fast algorithms for EAFE nonlinear discrete systems to solve a class of general time-dependent PNP equations.
To overcome the drawbacks such as slow convergence speed of the Gummel iteration, a geometric FAS algorithm is proposed to improve the solving efficiency. Then, an algebraic FAS algorithm is designed to improve the limitations such as the requirement of nested coarse grids of the geometric FAS algorithm.
The numerical experiments including practical PNP equations in ion channel show that algebraic FAS can improve the efficiency of decoupling external iterations and expand the size of domain can solved compared with the classical Gummel iteration. Finally, the improved algorithms based on the acceleration technique and adaptive method to solve the problem of excessive coarse grid iterations in FAS algorithms when convective intensity is high. Numerical experiments have shown that the improved algorithms converge faster and significantly enhance the adaptability of convective intensity. These fast algorithms presented in this paper is easy to extend to other FE nonlinear discrete systems, such as classical FE and SUPG method etc. It is promising to apply them to more complex PNP models, such as PNP equations in
more complex ion channels and semiconductor devices.

\noindent
{\bf Acknowledgement } S. Shu was supported by the China NSF (NSFC 12371373). Y. Yang was supported by the China NSF (NSFC 12161026).

\bibliographystyle{unsrt}  
\bibliography{artpnp}
\end{document}